\documentclass[11pt,twoside]{article}
\usepackage{latexsym}
\usepackage{amssymb,amsbsy,amsmath,amsfonts,amssymb,amscd}
\usepackage{graphicx,color}
\usepackage{float,url}
\usepackage{cancel}
\usepackage[colorlinks,linkcolor=red,citecolor=blue]{hyperref}

\newcommand  \ind[1]  {   {1\hspace{-1.2mm}{\rm I}}_{\{#1\} }    }
\newcommand{\R}{\mathbb{R}}
\newcommand {\e}  {\varepsilon}
\newcommand {\sg} {\sigma}
\newcommand {\vp} {\varphi}
\newcommand {\dv}  { {\rm div} }
\newcommand {\caL} { {\mathcal L} }
\newcommand {\cT} { {\mathcal T} }
\newcommand {\p}   {\partial}
\newcommand{\dis}{\displaystyle}

\newcommand{\beq}{\begin{equation}}
\newcommand{\eeq}{\end{equation}}
\newcommand{\bea} {\begin{array}{rl}}
\newcommand{\eea} {\end{array}}
\newcommand{\bepa}{\left\{ \begin{array}{l}}
\newcommand{\eepa} {\end{array}\right.}

\newcommand{\intt}{\int \hskip-8pt\int}
\newcommand{\intttt}{\int \hskip-8pt\int  \hskip-8pt\int \hskip-8pt\int}

\newcommand{\vip}{\vskip0.15cm}

\newtheorem{theorem}{Theorem}%[section]

\newtheorem{remark}[theorem]{Remark}
\newtheorem{proposition}[theorem]{Proposition}

\title{Transport distances for PDEs: the coupling method}

\author{
Nicolas Fournier\thanks{Sorbonne Universit\'{e}, CNRS, Laboratoire de Probabilit\'e, Statistique et Mod\'elisation, F-75005 Paris, France. Email: Nicolas.Fournier@sorbonne-universite.fr}
\and \addtocounter{footnote}{5}
Beno\^\i t Perthame\thanks{Sorbonne Universit\'{e}, CNRS, Universit\'{e} de Paris, Inria, Laboratoire Jacques-Louis Lions, F-75005 Paris, France. 
Email: Benoit.Perthame@sorbonne-universite.fr. 
B.P. has received funding from the European Research Council (ERC) under the European Union's Horizon 2020 research and innovation programme (grant agreement No 740623). }
}
\date{\today}

\begin{document}
\maketitle
\pagestyle{plain}
\pagenumbering{arabic}

\begin{abstract} 
We informally review a few PDEs for which some transport cost between pairs of solutions, 
possibly with some judicious cost function, decays: heat equation, Fokker-Planck equation,
heat equation with varying coefficients, fractional heat equation with varying coefficients,
homogeneous Boltzmann equation for Maxwell molecules, and some nonlinear integro-differential equations
arising in neurosciences.
We always use the same method, that consists in building a coupling between two solutions. This  means that we double the variables and solve, globally in time,  
a well-chosen PDE posed on the Euclidian square of the physical space.  
Finally, although the above method fails, we recall a simple idea to treat the case of the porous media equation. We also introduce another method based on the dual Monge-Kantorovich problem.
\end{abstract} 
%\vskip .7cm

\noindent{\makebox[1in]\hrulefill}\newline
2010 \textit{Mathematics Subject Classification:} 35A05; 35K55; 60J99; 28A33 % 35K55, 35B25, 76D27, 92C50.
\newline\textit{Keywords and phrases:} Transport distances; Monge-Kantorovich distance; Coupling; 
Fokker-Planck equation; Fractional Laplacian; Homogeneous Boltzmann equation; Integro-differential equations;
Porous media equation.\footnote{\textit{Acknowledgments:} We warmly thank the referees for their numerous fruitful comments.}

\tableofcontents

%
%%%%%%%%%%%%%%%%%%%%%%%%%%%%%%%%%%%%%%%%%%%%
\section*{Introduction}
%\label{sec:intro}
%-------------------------------------------
%%%%%%%%%%%%%%%%%%%%%%%%%%%%%%%%%%%%%%%%%%%%

It is usual to study the well-posedness, stability and large-time behavior of stochastic processes (e.g. solutions
to Stochastic Differential Equations) by using coupling methods: we consider two such processes, with different initial conditions,
driven by {\it suitably correlated} randomness, and we measure the 
$\varrho$-transport cost $\cT_{\varrho}$ between their distributions.

\vip

We work in $\R^d$ and we always assume that the cost function $\varrho:\R^d\times\R^d\mapsto \R$ satisfies
$\varrho (x,x)=0$ and $\varrho (x,y)=\varrho (y,x)>0$ for $x \neq y$. We recall that 
for two probability densities $u_1, u_2$ on $\R^d$,
\beq\bepa
\cT_{\varrho}(u_1, u_2) = \dis \inf_{v \in K(u_1,u_2)} \intt \varrho(x,y) v(x,y) dx \, dy, \\[10pt]
K(u_1,u_2)= \{v: \R^d\times\R^d\mapsto \R_+ \quad \hbox{such that} \quad \int v(x,y) dy=u_1(x), \;  \int v(x,y) dx =u_2(y) \}.
\eepa
\label{eq:wasserstein} \eeq
When for some $p\geq 1$,
$$
\varrho_p (x,y) =  \frac{|x-y|^p}{p},
$$
and we put $\cT_p=\cT_{\varrho_p}$. 
The distances $d_p=\cT_p^{1/p}$ are also called Wasserstein distances, and one refers
to the Monge-Kantorovich distance when $p=1$. Recent accounts about the theory can be found in the books \cite{AGSbook,VTOT, car}.

\vip

The probabilistic coupling method consists in finding, for two solutions $u_1$ and $u_2$ of a given 
Partial Differential Equation (PDE in short), two coupled stochastic processes $X_1$ and $X_2$, with time-marginals
$u_1$ and $u_2$, so that $X_1$ and $X_2$ remain as close as possible. One then controls
$\cT_\varrho(u_1(t),u_2(t))$ by $\mathbb{E}[\varrho(X_1(t),X_2(t))]$.
The goal of the present survey paper is to describe, in an informal way, this method, 
using only arguments based PDEs. This is what we call the {\it coupling method}: 
for each problem, we introduce a PDE, with doubled variables, describing the time-evolution of density 
$v(x,y,t)$ of the law at time $t\geq 0$ of two underlying coupled processes, 
in such a way that its two marginals
solve the original equation, and such that $v(x,y,t)$ is as concentrated as possible near the diagonal $x=y$.
This is possible only when the underlying pair of stochastic coupled processes is Markov. 
We then study the quantity
$\iint \varrho(x,y)v(x,y,t)$, which controls the transport cost $\cT_\varrho$ between two solutions.
This can be seen as a translation, in terms of PDEs, of the probabilistic coupling method.
For example, our proof of Tanaka's theorem  \cite{Tanaka_b} for the Boltzmann equation
really relies on the same main arguments as Tanaka (who is using some Poisson-driven stochastic differential
equations, but we avoid introducing any stochastic process).

\vip

The difficulties and novelties rely on the choice of the cost function and on the choice of coupling 
between two solutions by solving
a well-chosen PDE posed on the Euclidian square of the physical space, $\R^{2d}$ in general. Each time, we try to
emphasize the main technical difficulties that would allow one to justify the computations.

\vip

Of course, any way to produce some non-expansion estimates along solutions of some PDEs 
for some transport cost relies
on some coupling. However, for example, the deterministic methods in Carrillo \cite{car} for the heat equation or
Villani \cite[Section 7.5.6]{VTOT}, see also \cite{BoCaBoltz}, for Tanaka's theorem, 
are really different in spirit. In \cite{car}, the main tool is that solutions to the heat equation 
can be represented by a convolution formula. In \cite[Section 7.5.6]{VTOT} and \cite{BoCaBoltz},
everything relies on the contractive property of the gain operator and on the Duhamel principle.

\vip

For example, considering the Brownian motion leads to the heat equation.
We first give a simple proof that the heat equation is non-expansive (weak contraction) for 
any smooth cost function of the form $\varrho(x,y)=r(|x-y|)$. This is standard but the PDE literature 
seems to ignore this simple approach. 
The method can be extended to various cases. The Fokker-Planck equation is the simplest extension.
The case of the heat equation with variable coefficients, of the form 
$\p_t u - \Delta (a(x) u ) =0$, is more involved:
in one dimension, the Monge-Kantorovich distance $\cT_1$ plays a central role and is always non-expansive (under technical conditions);
we illustrate the general structure in higher dimension and show that if the cost function $\varrho$ 
satisfies some elliptic PDE, which does seem to enter a class with generic existence results, then $\cT_\varrho$ is non-expansive along solutions. 
The method also applies to some jump processes: fractional heat equation with variable coefficients in dimension
one, scattering equations, kinetic scattering equations, Boltzmann equation for Maxwell molecules.

\vip

For the porous media equation, the situation is more intricate and the above method does not seem to apply.
However, we recall  from \cite {BoCaPM} another, somehow related and rather simple, path to treat this equation. 

\vip

Concerning piecewise deterministic jump processes and (inhomogeneous) kinetic scattering equations, we present a new result showing that  some transport costs are non-expansive.

\vip

Finally, concerning jump processes, in particular those related to the discretized heat equation, we present another approach, based on the dual formulation of the transport costs.

\vip

The first example of use of the coupling method, up to our knowledge, can be traced back to Dobrushin \cite{Dob}, where the 
Vlasov equation is derived as mean-field limit of a deterministic system of interacting particles,
making use of some transport cost.
No PDE is written for the coupling in \cite{Dob}, because everything may be written in terms of 
characteristics. See \cite[Section 3]{GoMoPa} for a PDE analogue of Dobrushin's argument. 
In the same spirit, the Euler equation is derived from a deterministic system of interacting vortices 
in Marchioro-Pulvirenti
\cite[Section 5.3]{MaPu}, using also a coupling argument. See \cite{Ha} for a result with the strong 
transport distance $d_\infty$ but for a problem slightly less singular than the Euler equation.

\vip

Recently, the topic of transport costs has developed quickly for PDEs and  integro-differential equations  (IDEs) after new  understanding of optimal transportation and of the Brenier-Kantorovich map by~\cite{Brenier_polar,BeBrComp}. There are several approaches to use the transport costs in PDEs.  A geometrical approach based on gradient flow structures has been introduced in~\cite{Otto_CPDE2001} and extended in \cite{CMV2006, BGG_2013}, in particular for the porous media equation with $\cT_2$, for interacting particle systems and for granular flows. Also, many results on PDEs have been derived from the variational discretization 
algorithm named JKO after~\cite{JKO_98}. See the book~\cite{V_OT2009} for a complete presentation of these results.  See also \cite{Vas} who showed that contractivity for some other transport cost may fail for the
porous media equation. 
Let us also mention that the special structure associated with dimension 1  has been used to prove strict contraction for the porous media equation~\cite{CDT_07} for  the cost $\cT_2$, and to treat other equations as scalar conservation laws~\cite{BBL2005}, the Keller-Segel system~\cite{Cal_Car2012}  or the granular 
media equation \cite{LiTo}. See also the general reference \cite{CaTo} for nonlinear diffusion equations
and \cite{BlCaCa} for the multi-dimensional Keller-Segel system.
Methods based on optimal transportation have also been recently used to treat singular 
{\it congestion} (incompressible) equations arising in crowd modeling, see for instance \cite{BrCaSa, MRS2, CrKimYao}.

\vip

Most of the recent papers using some transport cost for PDEs have been using the gradient flow structure which is closely related to a variational formulation of the fluxes.
Here, with several examples of conservative equations, which do not necessarily have a gradient flow structure, we control the transport cost using the coupling method. We often borrow our examples from the stochastic processes which represent the PDEs thanks to their Kolmogorov equation. The cases of variable coefficients are particularly interesting because they often require some special
choice of the cost function.

\vip 

We organize our examples as follows. We begin  with three simple examples: heat equation, Fokker-Planck equation, and a class
of nonlinear transport equations. We show directly that the transport costs are non-expansive along these equations. 
Then we turn, in Section \ref{sec:hvc}, to the heat equation with variable coefficients. In Section~\ref{sec:integral} we  consider some
IDEs: scattering equations, including kinetic scattering and inhomogeneous fractional heat equation. 
Zero-th order terms, describing absorption and re-emission, as they appear in models of neural networks, can also be treated by adapting the method; this is explained in Section~\ref{sec:inhomogeneous}. The famous Tanaka theorem for the homogeneous Boltzmann equation can be included in our framework and this is done in Section~\ref{sec:Beq}. We treat the porous media equation in Section~\ref{sec:PME}. Finally, we exemplify in Section~\ref{sec:duality} how the same results can be proved using the dual formulation of the transport costs.

%%%%%%%%%%%%%%%%%%%%%%%%%%%%%%%%%%%%%%%%%%%%
\section{Heat, Fokker-Planck and transport equations}
\label{sec:HFP}
%-------------------------------------------
%%%%%%%%%%%%%%%%%%%%%%%%%%%%%%%%%%%%%%%%%%%%

In order to explain the coupling method in a very simple, but still relevant, framework, we  begin with 
the heat equation. Then we turn to drift and transport terms.

%------------------------------------------------------------------------
\subsection{Heat equation}
\label{sec:heat}
%-------------------------------------------------------------------------

Here is the well-known result, see e.g. \cite{V_OT2009}, we want to quickly recall.
\begin{theorem} 
Consider any increasing function $r:[0,\infty)\mapsto[0,\infty)$ such that the cost function
$\varrho:\R^d\times\R^d\mapsto \R_+$ defined by $\varrho(x,y)=r(|x-y|)$ is of class $C^2$.
Consider two probability densities $u^0_1, u^0_2$ on $\R^d$, and the corresponding solutions
$u_1,u_2$ to the heat equation 
\beq
\p_t u - \Delta u =0, \qquad x \in \R^d, \; t \geq 0.
\label{heat} 
\eeq
For any $t\geq 0$,  one has
$$
\cT_{\varrho}(u_1(t), u_2(t)) \leq \cT_\varrho(u_1^0, u_2^0).
$$
\label{th:H}
\end{theorem}

\begin{proof} We consider an initial density $v^0:\R^d\times\R^d\mapsto \R_+$ with marginals 
$u^0_1$ and $u^0_2$, that is such that $v^0 \in K(u_1^0,u_2^0)$.
We next consider the solution $v(x,y,t)$ to the degenerate heat equation
\beq
\frac{\p  v}{\p t} -\Delta_x v -\Delta_y v - 2 \nabla_x\cdot \nabla_y v =0, \quad x, y \in \R^d, \; t \geq 0
\label{heat:v}
\eeq
starting from $v^0$. Observe that \eqref{heat:v} may be rewritten as
\beq\label{hhh}
\frac{\p  v}{\p t} - (\nabla_x+\nabla_y)\cdot (\nabla_x+\nabla_y) v =0, \quad x, y \in \R^d, \; t \geq 0.
\eeq
Clearly, it holds that $v(x,y,t) \geq 0$, because of the
non-negativity of the operator in \eqref{heat:v}, which can be written in the variables 
$(x+y, x-y)$ as $-\Delta_{x+y}$.

\vip

We then define the marginals
$$
v_1(x,t) =\int v (x,y,t) dy, \qquad v_2(y,t) =\int v (x,y,t) dx
$$
and show that $v_1=u_1$ and $v_2=u_2$: for instance, integrating  \eqref{heat:v} with respect to $y$,
one finds
$$
\bepa
\frac{\p v_1(x,t) }{\p t} -\Delta_x v_1(x,t) =0, \quad x\in \R^d, \; t \geq 0,
\\[5pt]
v_1(x) = u^0_1(x), \quad x\in \R^d
\eepa
$$
and uniqueness of the solution of the heat equation gives us $v_1=u_1$.

Recalling \eqref{eq:wasserstein}, we conclude that 
$$
\cT_\varrho(u_1(t),u_2(t))\leq \intt \varrho(x,y)v(x,y,t)dxdy=\intt r(|x-y|)v(x,y,t)dxdy.
$$

Finally, we may also compute, using \eqref{hhh} and integrating by parts,
$$
\frac{d }{d t} \intt  r(|x-y|) v(x,y,t) dx dy = \! \intt v(x,y,t) \Big((\nabla_x+\nabla_y)\cdot (\nabla_x+\nabla_y)
[r(|x-y|)] \Big)dx dy =\!0.
$$

Therefore, for any initial data $v^0\in K(u_1^0,u_2^0)$,
$$
\cT_{\varrho} (u_1(t), u_2(t)) \leq \intt r(|x-y|)v^0(x,y)dxdy= \intt \varrho(x,y)v^0(x,y)dxdy
$$
and minimizing among such $v^0$ completes the proof. 
\end{proof} 

\vip

The only technical question is to justify the integration by parts, which is immediate if we assume 
enough moments initially (otherwise there is nothing to prove), at least 
when we restrict ourselves to power cost functions $\varrho_p(x,y)=|x-y|^p/p$ with $p\geq 2$. 
Notice that the well-posedness for~\eqref{heat:v}  follows from the observation that we actually
deal with $-\Delta_{x+y}$.
It is also possible, under some conditions, to treat the case of some non smooth cost functions, e.g. 
$\varrho_p$ for some $p\in [1,2)$: this issue is discussed in Section~\ref{sec:hvc}.

\vip

Remark that the solution to \eqref{heat:v} can be represented as $v(x,y,t)=\mathbb{E}[v^0(x+
\sqrt 2 B_t,y+\sqrt 2 B_t)]$,
for $(B_t)_{t\geq 0}$ a Brownian motion. Note also that in  \eqref{heat:v}, diffusion occurs only in
the $x+y$ variable, while $x-y$ remains invariant under the dynamics.

\vip

Another way to understand \eqref{heat:v} is to see that $v(t)$ is the law of the couple
$(X_0\!+\!\sqrt 2 B_t,Y_0\!+\!\sqrt 2B_t)$, where $(X_0,Y_0)$ is $v^0$-distributed and where we use the {\it same} 
Brownian motion $(B_t)_{t\geq 0}$ for both coordinates. 
If using e.g. independent Brownian motions $(B_t)_{t\geq 0}$ and $(C_t)_{t\geq 0}$, 
i.e. if considering the coupling $(X_0+\sqrt 2B_t,Y_0+\sqrt 2C_t)$, one would find $v$ with the good marginals but solving
$\partial_t v - (\Delta_x+\Delta_y)v=0$, leading for example to the far from optimal estimate 
$\cT_2(u_1(t),u_2(t))\leq \cT_2(u_1^0,u_2^0)+2d t$.

%------------------------------------------------------------------------
\subsection{Fokker-Planck equation}
%-------------------------------------------------------------------------

The coupling method can be extended to the Fokker-Planck equation, see \cite{BCG2012} for some 
more elaborate consequences. The result can be stated as follows

\begin{theorem} 
Consider some function $V:\R^d\times \R_+\mapsto \R^d$ such that, for some $\alpha\in\R$,
\beq
(V(x,t) - V(y,t))\cdot(x-y) \leq \alpha |x-y|^2, \qquad   x,  y \in \R^d, \; t \geq 0.
\label{FP_as} \eeq
Consider two probability densities $u^0_1, u^0_2$ on $\R^d$ and the corresponding solutions
$u_1,u_2$ to the Fokker-Planck equation
\beq
\p_t u - \Delta u + \dv(V(x,t) u)=0, \qquad x \in \R^d, \; t \geq 0.
\label{FP} \eeq
For any $t\geq 0$, any $p\geq 1$, one has
$$
\cT_p(u_1(t), u_2(t)) \leq \cT_p(u_1^0, u_2^0) \exp(\alpha p t).
$$
\label{th:FP}
\end{theorem}

This inequality is  well-known, see for example \cite{V_OT2009} \S 9.1.5 and the references therein. One can also find relations to several deep and recent  functional analysis tools. This goes far beyond our present purpose.

\vip

\begin{proof}
This is the same proof as for the heat equation, with longer expressions.
We consider any $v^0$, with marginals $u^0_1$ and $u^0_2$, and the solution
$v$ to the equation
\beq
\p_t v - \Delta_x v  -  \Delta_y v - 2 \nabla_x\cdot \nabla_y v + \dv_x(V(x,t) v) + \dv_y(V(y,t) v) =0, \qquad x,  y \in \R^d, \; t \geq 0
\label{FP:v} \eeq
starting from $v^0$.
One easily checks that $v(x,y,t)\geq 0$ and, integrating \eqref{FP:v} with respect to $y$,
that $v_1(x,t):=\int v(x,y,t)dy$ solves \eqref{FP} and starts from $u^0_1$,
whence $v_1=u_1$. The second marginal is treated similarly, and we conclude that
$\cT_p(u_1(t),u_2(t))\leq p^{-1}\int\hskip-4pt\int |x-y|^p v(x,y,t)dxdy$.
Finally, using the same computation as for the heat equation, with some additional terms, we see that
$$\bea
\dis \frac{d }{d t} \intt  \frac{|x-y|^p}{p} v(x,y,t) dx dy 
&= 0+\dis  \intt v(x,y,t) |x-y|^{p-2}  (x-y)\cdot (V(x,t)-V(y,t)) dx dy \\ [10pt]
&\leq \dis  \alpha \intt  |x-y|^p v(x,y,t) dx dy
\eea$$
by assumption \eqref{FP_as}.
The result follows using the Gronwall lemma
$$
\cT_p(u_1(t),u_2(t))\leq \intt \frac{|x-y|^p}{p} v(x,y,t)dxdy \leq \Big( \intt \frac{|x-y|^p}{p} v^0(x,y)dxdy\Big)e^{\alpha p t}
$$
and minimizing in $v^0$.
\end{proof}

\vip

Observe that one can treat in a similar way the case with variable diffusion coefficients,
of the form $\partial_t u-\sum_{i,j=1}^d \partial_{ij} (a_{ij}(x,t) u)+ \dv(V(x,t) u)=0$, under suitable
conditions on the nonnegative symmetric matrix $a$ and on $V$. For example, we will find that
$\cT_2(u_1(t),u_2(t))\leq \cT_2(u_1^0,u_2^0)\exp(\alpha t)$ as soon as 
$$
{\rm Tr} \Big([\sigma(x,s)-\sigma(y,s)][\sigma(x,s)-\sigma(y,s)]^*\Big) + (x-y)\cdot (V(x,t)-V(y,t))
\leq \alpha |x-y|^2,
$$
where $\sigma(x,t)$ is, for each 
$(x,t)\in \R^d\times[0,\infty)$, a matrix such that $\sigma(x,t)[\sigma(x,t)]^*=a(x,t)$.
\vip
To our knowledge, a divergence form equation
$\partial_t u-\sum_{i,j=1}^d \partial_{i} (a_{ij}(x,t) \partial_j u)=0$ does not enjoy particular
properties from this point of view.

%------------------------------------------------------------------------
\subsection{A nonlinear transport equation}
%-------------------------------------------------------------------------

We next consider a fully deterministic problem  which 
arises in several types of modeling, such as polymers, cell division, neuron networks, etc :
\beq
\p_t u + \dv [V(x,I(t)) u] =0, \qquad x \in \R^d, \; t \geq 0,
\label{NLTR}
\eeq
where the nonlinearity stems from the quantity $I(t)$ defined, with a given weight $\psi:\R^d\mapsto \R$,
\beq
I(t) = \int_{\R^d} \psi(x) u(x,t) dx.
\label{NLTR2}
\eeq
We again complement this equation with an  initial condition $u^0 \geq 0$ with mass $\int  u^0=1$.
%--------------------------------
\begin{theorem} \label{ettac}
Assume that $V:\R^d\times\R \mapsto \R^d$ and $\psi:\R^d\mapsto \R$ are of class $C^1$, 
and that for some $\alpha>0$
\beq
(x-y)\cdot\big(V(x,I)-V(y,I)\big) \leq - \alpha |x-y|^2, \qquad \forall x, y \in \R^d, \; I \geq 0.
\label{NLTR3}
\eeq
Setting $\langle x(t) \rangle = \int x u(x,t) dx$, we have
\[
\cT_2(u(t),\delta_{\langle x(t) \rangle})=\int \frac{ |x- \langle x(t) \rangle |^2}{2}  u(x,t)  dx  
\leq e^{-2\alpha t}  \int \frac{|x-\langle x(0) \rangle|^2} {2} u^0(x,t)dx = e^{-2\alpha t}\cT_2(u(0),\delta_{\langle x(0) \rangle}).
\]
Assume additionally that
$$
\beta= \| D_I V \|_\infty  \| D\psi \|_\infty < \alpha 
$$
and fix any initial point $X^0\in\R^d$.
Consider the solution $X$ to $X'(t)=-V(X(t),\psi(X(t))$ starting from $X^0$.
For all $t\geq 0$, one has
$$
\cT_2(u(t),\delta_{X(t)})= \int \frac{|x-X(t)|^2}{2} u(x,t) dx \leq  e^{2 (\beta-\alpha) t} \int \frac{|x-X^0|^2}{2} u^0(x) dx
=e^{2 (\beta-\alpha) t} \cT_2(u^0,\delta_{X^0}).
$$
\end{theorem}
%--------------------------------
It holds that $(\delta_{X(t)})_{t\geq 0}$ solves \eqref{NLTR} in a weak sense.
A more general result, involving any pair of solutions, can be found in \cite{V_OT2009}.
One could prove Theorem \ref{ettac} without using a PDE for the coupling, using only the characteristics, 
in the spirit of Dobrushin \cite{Dob}.

\vip

\begin{proof}  We consider two solutions $u_1$ and $u_2$ to \eqref{NLTR}, and denote by $I_1(t)$ and $I_2(t)$ 
the corresponding functions, see \eqref{NLTR2}. As we are interested in the case where one of the two solutions
is a Dirac mass (for each $t\geq 0$, $u_2(t)=\delta_{X(t)}$), we can consider the trivial coupling
$v(x,y,t)= u_1(x,t) u_2(y,t)$, which of course has the correct marginals, and satisfies
$$
\p_t v + \dv_x [V(x,I_1(t)) v] +  \dv_y [V(y,I_2(t)) v] =0.
$$
Therefore, we may compute
$$\bea
\dis\frac d {dt}\intt \frac{ |x-y|^2}{2}& u_1(x,t) u_2(y,t) dx dy 
= \dis \intt (x-y)\cdot\big(V(x,I_1(t))-V(y,I_2(t))\big) u_1(x,t) u_2(y,t)  dx dy
\\[10pt]
=&\dis  \intt (x-y)\cdot\big(V(x,I_1(t))-V(y,I_1(t))\big) u_1(x,t)u_2(y,t)  dx dy
\\[10pt]
&+\dis  \intt (x-y)\cdot\big(V(y,I_1(t)) ) - V(y,I_2(t))\big) u_1(x,t)u_2(y,t)  dx dy
\\[10pt]
\leq& -\dis \alpha \intt\! |x-y|^2 u_1(x,t)u_2(y,t) dx dy \\
&+\dis  \| D_I V \|_\infty | I_1(t)-I_2(t)| 
\Big(\intt |x-y|^2 u_1(x,t)u_2(y,t) dx dy \Big)^{1/2}.
\eea
$$
We first apply this in the case of single solution $u:=u_1=u_2$, whence $I_1=I_2$, and
we directly conclude by Gronwall's lemma that 
\[
\intt \frac{ |x-y|^2}{2} u(x,t) u(y,t) dx dy \leq e^{-2 \alpha t} \intt \frac{ |x-y|^2}{2} u^0(x,t) u^0(y,t) dx dy.
\]
This classically rewrites as
\[
\int |x- \langle x(t)\rangle |^2  u(x,t)  dx \leq e^{-2 \alpha t}  \int { |x-\langle x(0)\rangle|^2} u^0(x,t) dx.
\]
as desired. Next, when considering two solutions,  we notice that
$$
I_1(t)-I_2(t)  = \intt [\psi(x)- \psi(y)]u_1(x,t)u_2(y,t) dx dy,
$$
whence
$$
\big| I_1(t)-I_2(t) \big| \leq \| D\psi \|_\infty  \intt |x - y| u_1(x,t)u_2(y,t) dx dy 
\leq \| D\psi \|_\infty \Big(\intt |x-y|^2 u_1(x,t)u_2(y,t) dx dy \Big)^{1/2}.
$$
Therefore,
$$
 \frac d {dt}\dis \intt \frac{ |x-y|^2}{2} u_1(x,t)u_2(y,t) dx dy 
\leq (\beta-\alpha)\dis \intt |x-y|^2 u_1(x,t)u_2(y,t) dx dy.
$$
Applying this to the case where $u_2(t)=\delta_{X(t)}$ concludes the proof.
\end{proof} 

\vip

 For consistency with the other presentations in this section, we have written this result for an $L^1$ density 
$v$ with a finite second moment, but the extension to a probability measure is immediate.

%%%%%%%%%%%%%%%%%%%%%%%%%%%%%%%%%%%%%%%%%%%%
\subsection{Aggregation type equation}
\label{sec:agregation}
%-------------------------------------------
%%%%%%%%%%%%%%%%%%%%%%%%%%%%%%%%%%%%%%%%%%%%

Among the nonlinear transport equations let us also mention the aggregation equation 
\beq
\p_t u - \dv [ u \, (K*u)] =0, \qquad x \in \R^d, \; t \geq 0,
\label{agr1}
\eeq
where the convolution  term is defined by  $K*u(x,t) = \int_{\R^d} K(x-y) u(y,t) dy$. It is well-known, under various assumptions on $K$ and for several variants of the equation, that it  is contractive for transport costs and we refer to the papers~\cite{Otto_CPDE2001,CMV2006,CaTo, LiTo, GoMoPa}.

Here we assume that the smooth kernel satisfies, for some $\alpha \geq 0$,
\beq
K\in C^1(\R^d; \R^d), \qquad  (x-y)\cdot\big(K(x)-K(y) \big) \geq - \alpha {|x-y|^2} \qquad K(x)=-K(-x).
\label{agr:as1}
\eeq
We also assume that the $L^1(\R^d)$-initial conditions $u_i^0 \geq 0$ with  unit mass $\int  u_i^0=1$, for $i=1, \, 2$ satisfy
\beq
\int_{\R^d} |x|^2 u^0_i dx < \infty, \qquad \int_{\R^d} x u^0_i (x) dx  =0.
\label{agr:as2}
\eeq
These imply that $ \int_{\R^d} x u_i(x,t) dx  =0$ for all $t>0$.

%--------------------------------
\begin{theorem} \label{th:agregation}
Assume \eqref{agr:as1}. Consider two initial data $u_1^0$ and $u_2^0$ satisfying \eqref{agr:as2} 
and $u_1$, $u_2$
the corresponding solutions to \eqref{agr1}. One has for all $t \geq 0$,
$$
\cT_2(u_1(t),u_2(t) ) \leq e^{2\alpha t} \cT_2(u_1^0,u_2^0).
$$
\end{theorem}

Notice that, combining the coupling below with that in Section~\ref{sec:heat}, it is immediate to also treat the case with diffusion
$$
\p_t u - \Delta u - \dv [ u \, K*u] =0,
$$
Notice also that when $K=\nabla U$ for some potential $U(x)$, then the equation has a gradient flow structure, which is not the case here.

\begin{proof}  Again the simple coupling 
\beq
\p_t v(x,y,t) - \dv_x [ v \, K*u_1] - \dv_y [ v \, K*u_2] =0, \qquad x,\,y \in \R^d, \; t \geq 0
\label{agr3}
\eeq
is enough.  Indeed, we can verify immediately the properties that it preserves nonnegativity and the marginales satisfy the correct equations. Finally, one can compute
$$\bea 
\dis \frac{d}{dt} \intt &\dis\frac{|x-y|^2}{2}  v(x,y,t) dx dy 
\\ & =\dis  -\intt (x-y)\cdot \Big[ \int K(x-x') u_1(x') dx'-   \int K(y-y') u_2(y') dy' \Big] v(x,y,t) dx dy
\\[12pt]
&\dis = - \intttt  (x-y)\cdot\big[ K(x-x') -  K(y-y') \big] v(x',y',t) v(x,y,t) dx' dy' dx dy 
\\[12pt]
&\dis = - \frac 12  \intttt  (x-x'-y+y')\cdot\big[ K(x-x') - K(y-y') \big] v(x',y',t) v(x,y,t) dx' dy' dx dy 
\\[12pt]
&\dis \leq  \frac  \alpha 2  \intttt |x-x'-y+y' |^2 v(x',y',t) v(x,y,t) dx' dy' dx dy 
\\[12pt]
&\dis = \alpha \intt {|x-y|^2} v(x,y,t)  dx dy.
\eea .
$$
We have used a symmetry argument for $K$ and that $\int\!\!\int (x-y)v(x,y,t) dx dy=0$.
We complete the proof as in Section~\ref{sec:heat}. 
\end{proof} 

%%%%%%%%%%%%%%%%%%%%%%%%
%-------------------------------------------------------------
\section{Heat equation with variable coefficients}
\label{sec:hvc}
%------------------------------------------------------------

We consider the heat equation with variable coefficient.
This is much more intricate than the previous examples.
In 1 dimension, we use the $\cT_1$ distance and recover a result implicitly included in \cite{LeGall_1983}.
In higher dimension, we indicate a general way to construct cost functions. This leads
to a poorly explored degenerate elliptic PDE, see however \cite{Pier2010} and the references therein.

%
%%%%%%%%%%%%%%%%%%%%%%%%%%%%%%%%%%%%%%%%%%%%
\subsection{One-dimensional case}
\label{sec:d1}
%-------------------------------------------
%%%%%%%%%%%%%%%%%%%%%%%%%%%%%%%%%%%%%%%%%%%%

We consider some $a:\R\mapsto \R_+$ and the following heat equation.
\begin{equation}
\frac{\p u}{\p t} -\frac{\p^2}{\p x^2} [a(x) u ] =0, \quad x\in \R, \; t \geq 0.
\label{eq:heatcoef} 
\end{equation}

\begin{theorem} Assume that $d=1$ and that $a= \sg^2$ for some $\sg \in C^{1/2}(\R)$.
Consider two probability densities $u^0_1,u^0_2$ on $\R$ and the corresponding solutions
$u_1,u_2$ to \eqref{eq:heatcoef}. For all $t\geq 0$, one has
$$
\cT_{1} (u_1(t), u_2(t)) \leq \cT_{1} (u_1^0, u_2^0).
$$ 
\end{theorem} 

\begin{proof}  We give a proof for $\sg \in C^\alpha(\R)$, with $\alpha >1/2$, the remark below explains how to treat 
$\alpha= 1/2$. We consider any probability density $v^0(x,y)$ with marginals $u_1^0$ and $u_2^0$ and consider the coupling equation
\beq 
\p_t v - \partial_{xx}(\sg^2(x) v) -\partial_{yy} (\sg^2(y) v) -2 \partial_x  \partial_y [ \sg(x)\sg(y) \; v] =0, 
\qquad x,  y \in \R, \; t \geq 0
\label{eq:join}
\eeq
starting from $v^0$. This equation preserves non-negativity. 
A simple way to see this is the following computation:
multiplying \eqref{eq:join} by $-v_-$, integrating on $\R^2$ and using some integrations by parts, 
one can check that 
\begin{align*}
&\frac 12 \frac{d}{dt} \intt v_-^2(x,y,t) dx \, dy \\
=& - \intt  \Big[ | \sg(x) \p_x v_-(x,y,t)|^2- 2 \sg(x) \sg(y) \p_x v_-(x,y,t)   \p_y v_-(x,y,t) +
| \sg(y) \p_y v_-(x,y,t)|^2 \Big]dx \, dy
\\
&+  \frac 12  \intt v_-^2(x,y,t)  [ \partial_{xx}(\sg^2(x))  + \partial_{yy}(\sg^2(y)) - 2 \p_x \sg(x) \p_y \sg(y) ] dx \, dy\\
\leq & \frac 12  \intt v_-^2(x,y,t)  [ \partial_{xx}(\sg^2(x))  + \partial_{yy}(\sg^2(y)) - 2 \p_x \sg(x) \p_y \sg(y) ] dx \, dy.
\end{align*}
Since $\int\hskip-4pt\int v_-^2(x,y,0) dx \, dy =0$, the result follows from the Gronwall lemma
if $\sg$ is smooth. Otherwise, one can work by approximation.
\vip

Integrating \eqref{eq:join} with respect to $y$, we see
that $v_1(x,t):=\int v(x,y,t)dy$ solves \eqref{eq:heatcoef} and starts from $u^0_1$,
whence $v_1=u_1$. The second marginal is treated similarly, and we conclude that
$\cT_1(u_1(t),u_2(t))\leq \int\hskip-4pt\int |x-y| v(x,y,t)dxdy$.
Because of its singularity, we need to regularize the absolute value as a $W^{2,\infty}$ function and define
$$
\omega_\e(r)=
\bepa
\frac{r^2}{2 \e} \qquad  \quad \; \hbox{ for } r\leq \e,
\\[5pt]
r - \frac \e 2 \qquad  \hbox{ for } r \geq \e .
\eepa
$$
Using the H\"older constant $C_\sg$ of $\sg(\cdot)$, we see that
$$\bea
\dis\frac{d}{dt}  \intt \omega_\e(|x-y|) v(x,y,t) dx dy& = \dis \intt  v(x,y,t) \;  \omega''_\e(|x-y|) 
\,  [\sg(x) - \sg(y)]^2 dxdy
\\[10pt]
&\leq C_\sg^2 \dis  \intt  v(x,y,t)  \frac{\ind{|x-y|\leq \e}}\e |x-y|^{2\alpha}dxdy 
\\[10pt]
&\leq C_\sg^2 \e^{2\alpha-1},
\eea $$
because $v(t)$ is a probability measure. Since now $2\alpha-1>0$, we may let $\e\to 0$ and we find that
$$
\cT_1(u_1(t),u_2(t))\leq \intt |x-y| v(x,y,t) dx dy \leq \intt |x-y| v^0(x,y) dx dy.
$$
We conclude, as usual, by minimizing in $v^0$.
\end{proof} 

\begin{remark} The condition $\sg \in C^{1/2}(\R)$ is enough. To treat this exponent,
a better construction of the regularization is required, using the so-called {\it Yamada function}:
\[
\omega_\e(r)= 0 \; \hbox{ for } r\leq \e^{3/2}, \qquad \omega_\e''(r)= \frac{2}{r |\ln(\e)|} \; \hbox{ for } \e^{3/2} \leq r\leq \e, \qquad 
\omega_\e'(r)=1 \;  \hbox{ for } r \geq \e .
\]
\end{remark}

There are other technical issues here. For example, the well-posedness of \eqref{eq:join},
which is necessary to identify the marginals of the solution $v$ to the coupling equation, is not so easy.
A possible direction is to use results established in \cite{fig}, in the spirit of \cite{dpl}.

%-------------------------------------------------------------
\subsection{A general construction of the weight}
\label{sec:genweight}
%------------------------------------------------------------

In order to unravel the algebraic structure behind the choice of the weight  $\varrho$, 
we now consider the general case of  dimension $d$.
We assume that $a:\R^d\mapsto {\mathcal M}_{d\times d}(\R)$ is everywhere
symmetric and nonnegative, of the form
\beq
a_{ij}(x) = \sum_{k=1}^K  \sg_{ik}(x) \sg_{jk} (x),
\eeq
for some $\sigma:\R^d\mapsto {\mathcal M}_{d\times K}(\R)$, and we consider the heat equation 
\beq 
\frac{\p u}{\p t} - \sum_{i,j=1}^d \frac{\p^2}{\p x_i \p x_j} [a_{ij}(x) u ] =0, \quad x\in \R^d, \; t \geq 0,
\label{eq:heatcoef_d} 
\eeq 
completed with an initial probability density $u^0$ on $\R^d$.

\begin{proposition} 
Assume that $\sigma$ is regular enough and 
consider two probability densities $u^0_1,u^0_2$ on $\R^d$ and the corresponding solutions
$u_1,u_2$ to \eqref{eq:heatcoef_d}. For all $t\geq 0$, one has 
$$
\cT_{\varrho} (u_1(t), u_2(t)) \leq \cT_{\varrho}(u_1^0, u_2^0) ,
$$
for any smooth cost $\varrho:\R^d\mapsto\R_+$ satisfying
\begin{equation}\label{eqp}
\sum_{i,j=1}^d  a_{ij}(x) \frac{\p^2  \varrho(x,y) }{\p x_i \p x_j} + \sum_{i,j=1}^d  a_{ij}(y) 
\frac{\p^2  \varrho(x,y) }{\p y_i \p y_j}  
 + 2 \sum_{i,j=1}^d   \sum_{k=1}^K \sg_{ik}(x) \sg_{jk} (y)  
\frac{\p^2 \varrho(x,y)}  { \p x_i \p y_j} \leq  0, \quad x,y \in \R^d.
\end{equation}
\end{proposition} 

When $a$ is constant, we recover that any $C^2$ cost function of the form $\varrho(x,y)=r(|x-y|)$ works. 
In dimension $1$, $\varrho(x,y)=|x-y|$ is indeed a (weak) solution to \eqref{eqp}. We do not know of a theory to solve \eqref{eqp}, in dimension $d\geq 2$, for a general coefficient $a$, so that we do not know if this result is useful. Notice that equation~\eqref{eqp} should be completed by the boundary value $ \varrho(x,x)=0$ with some growth condition to mimic $|x-y|^p$.
\vip

\begin{proof} We consider any probability density $v^0(x,y)$ with marginals $u_1^0$ and $u_2^0$ 
and consider the coupling equation
$$
\frac{\p  v}{\p t} -  \sum_{i,j=1}^d \frac{\p^2}{\p x_i \p x_j} [a_{ij}(x) v ]  -  \sum_{i,j=1}^d \frac{\p^2}{\p y_i \p y_j} [a_{ij}(y) v ] - 2 \sum_{i,j=1}^d \sum_{k=1}^K \frac{\p^2}{\p x_i \p y_j} [\sg_{ik}(x) \sg_{jk} (y)  v ]=0 
$$
starting from $v^0$. We show as usual that $\int v(x,y,t)dy=u_1(x,t)$ and that $\int v(x,y,t)dx=u_2(y,t)$.
Moreover, we have $v(x,y,t) \geq 0$:
we multiply the coupling equation by $-v_-$ and integrate, finding
\[ 
\frac 12 \frac {d}{dt} \intt  v_-^2(x,y,t) dx dy  = - \sum_{k=1}^K I_k  - J  ,
\]
with
\begin{align*}
I_k = \intt   \sum_{i,j=1}^d \Big[ \sg_{ik}(x)\frac{\p v_-(x,y,t)}{\p x_i}  \sg_{jk}(x) \frac{\p v_-(x,y,t)}{\p x_j } +  \sg_{ik}(y)\frac{\p v_-(x,y,t)}{\p y_i}  \sg_{jk}(y) \frac{\p v_-(x,y,t)}{\p y_j }\\
+ 2 \sg_{ik}(x) \sg_{jk}(y) \frac{\p v_-(x,y,t)}{\p x_i} \frac{\p v_-(x,y,t)}{\p y_j }\Big]dxdy
\end{align*}
which can also be written
\[
I_k = \intt \Big|   \sum_{i=1}^d \sg_{ik}(x)\frac{\p v_-(x,y,t)}{\p x_i} + \sum_{i=1}^d \sg_{ik}(y)\frac{\p v_-(x,y,t)}{\p y_i} \Big|^2 
dxdy\geq 0 .
\]
The other term is
\begin{align*}
J = \intt  \sum_{i,j=1}^d   \Big[  \frac{\p v_-(x,y,t)}{\p x_i}  \frac{\p a_{ij}(x)}{\p x_j}  v_-(x,y,t)+  \frac{\p v_-(x,y,t)}{\p y_i}  \frac{\p a_{ij}(y)}{\p y_j}  v_-(x,y,t) + \\ 
2  \frac{\p v_-(x,y,t)}{\p x_i}v_-(x,y,t) \frac{\p }{\p y_j}  \sum_{k=1}^K \sg_{ik}(x) \sg_{jk} (y)   \Big]dxdy ,
\end{align*}
which can also be written after integration by parts 
\[
J =  - \frac 12 \intt   (v_-(x,y,t))^2 \sum_{i,j=1}^d   \Big[   \frac{\p^2 a_{ij}(x)}{\p x_i\p  x_j}   +    \frac{\p ^2a_{ij}(y)}{\p y_i \p y_j} + 
2 \frac{\p^2 }{\p x_i \p y_j}  \sum_{k=1}^K \sg_{ik}(x) \sg_{jk} (y)   \Big]dxdy .
\]
Assuming that the entries $\sg_{ik}$ are bounded with two bounded derivatives, 
we conclude by Gronwall's lemma that $v_- \equiv 0$, since we initially have 
$\int\hskip-4pt\int v_-^2(x,y,0) dx \, dy =0$.

\vip

Recalling \eqref{eq:wasserstein}, we conclude that $\cT_{\varrho}(u_1(t),u_2(t))\leq \int\hskip-4pt\int \varrho(x,y) v(x,y,t) dx dy$.
Since finally
$$\bea
&\dis \frac{d }{d t}  \dis \intt  \varrho(x,y) v(x,y,t) dx dy 
\\[5pt]
=& \dis \intt v(x,y,t)  \sum_{i,j=1}^d  \Big[ a_{ij}(x) \frac{\p^2  \varrho(x,y) }{\p x_i \p x_j} +    a_{ij}(y) \frac{\p^2  \varrho(x,y) }{\p y_i \p y_j}  
 + 2 \ \sum_{k=1}^K \sg_{ik}(x) \sg_{jk} (y)  \frac{\p^2 \varrho(x,y)}  { \p x_i \p y_j} \Big] dxdy\leq 0
\eea$$
by assumption, we conclude as usual.
\end{proof} 

\vip

%\begin{remark} The coupled equation, under the form \eqref{eq:join}, 
%is the Kolmogorov equation for the coupled process
%\[ 
%dX_t  = \sg(X_t) dB_t, \quad dY_t  \;= \sg(Y_t ) \,dB_t .
%\]
%Hence the existence of a weight $\varrho$ is related to the 
%strong uniqueness of the process.
%\end{remark}

We leave open the question to formalize this approach rigorously, in particular for degenerate 
coefficients $\sigma$, and to build other examples where one can prove the existence of a weight $\varrho$.

%%%%%%%%%%%%%%%%%%%%%%%%
%-------------------------------------------------------------
\section{Scattering and integral kernels}
\label{sec:integral}
%------------------------------------------------------------

We now turn to equations that describe the probability law of various jump processes. These are well-known 
results except the case of kinetic scattering in Subsection~\ref{sec:kscattering} which seems to be new.

%-------------------------------------------------------------
\subsection{Simple scattering}
\label{sec:jump}
%------------------------------------------------------------

For $x\in \R^d$, we parameterize the pre-jump location $X= \Phi (x,h)$  by $h \in \R^d$, distributed according to a 
bounded measure $\mu$. We assume that for all fixed $h\in \R^d$, 
\beq
x\mapsto X= \Phi  (x,h) \hbox{ is invertible on $\R^d$ and }  D_x  \Phi (x,h) \hbox{ is an invertible matrix},
\label{jumpas1}\eeq
and we use the notation $X \mapsto x=  \Phi ^{-1}(X,h)$ for the inverse in $x$ (with $h$ fixed).

\vip

We consider the scattering problem
\beq
\p_t u(x,t)=\dis \int \Big[u( \Phi (x,h),t) \det (D_x  \Phi (x,h))  -u(x,t)\Big] d\mu(h),
\label{eq:jump}
\eeq
with initial condition $u^0$, a probability density on $\R^d$.
Actually, this equation is to be understood in the weak sense: integrating the right hand side 
against a test function $\vp(x)$, we see that
$$
\dis \intt \vp(x) \Big[u( \Phi (x,h),t) \det (D_x  \Phi (x,h))  -u(x,t)\Big] d\mu(h) = \dis \intt u(X,t) 
[\vp ( \Phi ^{-1}(X ,h) )- \vp(X)  ]  dX d\mu(h),
$$
which shows that the determinant $\det (D_x  \Phi (x,h)) $ is only used informally. 
We briefly prove the following result, which is classical, see for instance \cite{BalJ2018}.

\begin{theorem} Assume \eqref{jumpas1}, fix $p\in [1, \infty)$ and suppose there 
is $\delta \in \R$ such  that for all $X, Y \in\R^d$,
\beq
 \int |  \Phi ^{-1}(X,h)- \Phi^{-1}(Y,h) |^p d\mu(h) \leq K L | X - Y |^p, \quad \hbox{where} \quad K=\mu(\R^d).
\label{jumpas3}
\eeq 
Consider two probability densities $u^0_1,u^0_2$ on $\R^d$ and the corresponding solutions
$u_1,u_2$ to \eqref{eq:jump}. For all $t\geq 0$, one has 
$$
\cT_p (u_1(t), u_2(t)) \leq e^{K(L -1)t} \cT_p(u_1^0, u_2^0)  ,
$$
\end{theorem}
%--------------------------------------------
%--------------------------------------------

The homogeneous scattering corresponds to $ \Phi (x,h) =x+h$ and obviously fulfills the above assumptions.

\vip

\begin{proof} 
For a probability density $v^0$ on $\R^d\times\R^d$ with marginals $u^0_1$ and $u^0_2$,
we consider the solution $v$ of the coupled equation built in such a way that  the jumps parameter $h$ is common to the two variables. Namely, we choose
\beq
\p_t v(x,y,t)= \dis \int [v( \Phi (x,h),  \Phi (y,h),t) \det (D_x  \Phi (x,h)) \det (D_x  \Phi (y,h))  -v(x,y,t) ] d\mu(h),
\label{eq:jumpv}
\eeq
starting from $v^0$. We clearly have $v \geq 0$, and integrating in $y$ and using the change of variable
$y \mapsto  \Phi (y,h)$, we find that $v_1(x,t ) = \int v(x,y,t) dy$ satisfies \eqref{eq:jump}.
Since it starts from $u^0_1$, we conclude that $v_1=u_1$. The second marginal is treated similarly, and we conclude
as usual that $\cT_p(u_1(t),u_2(t))\leq p^{-1}\int\hskip-4pt\int |x-y|^pv(x,y,t)dxdy$.
Next, we compute, using~\eqref{eq:jumpv}:
$$\bea
\dis \frac d {dt} \intt & \dis |x-y|^p v(x,y,t)dx dy  + K\intt  |x-y|^p v(x,y,t) dxdy
\\
 &= \dis \intt \hskip-8pt \int |x-y|^p  v( \Phi (x,h),  \Phi (y,h),t) \det (D_x  \Phi (x,h)) \det (D_y  \Phi (y,h)) d\mu(h) dx dy
\\
&= \dis\intt\hskip-8pt \int | \Phi ^{-1}(X,h) - \Phi^{-1}(Y,h)|^p v(X,Y,t) d\mu(h)dX dY.
\eea $$
We used the changes of variables $X=  \Phi (x,h)$ and $Y= \Phi (y,h)$ (with $h$ fixed). Recalling~\eqref{jumpas3}, 
we conclude that 
$$
\frac d {dt} \intt \dis |x-y|^p v(x,y,t)dx dy \leq  K(L -1) \intt \dis |x-y|^p v(x,y,t)dx dy.
$$
Using the Gronwall lemma, we thus find that
$$
\cT_p(u_1(t),u_2(t)) \leq p^{-1}\intt |x-y|^pv(x,y,t)dxdy \leq  p^{-1}e^{K(L-1)t}\intt |x-y|^pv^0(x,y)dxdy
$$
and we conclude as usual, minimizing in $v^0$.
\end{proof} 

\vip

The most general scattering equation reads
\beq
\p_t u(x,t)=\dis \int \big[ \pi (x, x_ *) u( x_*)   -  \pi (x_ *,x) u(x,t)\big] dx_*,
\label{eq:jumpGen}
\eeq
and equation \eqref{eq:jump} corresponds to the homogeneous cases when $\int \pi(x_*,x)dx_*=1$, and the above method can easily be adapted. For the inhomogeneous case, see Section~\ref{sec:inhomogeneous}.

%-------------------------------------------------------------
\subsection{Kinetic scattering}
\label{sec:kscattering}
%------------------------------------------------------------

We next consider some kinetic scattering models, that means we work in the phase space.
We consider some finite measure $\mu$ on $\R^d$, some application $V:\R^d\mapsto\R^d$ such that, for all $h\in \R^d$,
\beq
v\mapsto V=  \Phi (v,h) \hbox{ is invertible and } D_v \Phi (v,h) \hbox{ is an invertible matrix},
\label{jumpas11}\eeq
and the kinetic scattering equation
\beq
\p_t f(x,v,t) +v.\nabla_x f = \int   [f(x, \Phi (v,h) ,t) \det (D_v  \Phi (v,h)) -f(x,v,t)]d\mu(h)
\label{KScattering}
\eeq
completed with an initial data $f^0(x,v) \geq 0$ with $\int f^0 dx dv =1$.

\begin{theorem} Assume \eqref{jumpas11}.
Set $K=\mu(\R^d)$ and suppose that for 
some $L\in \R_+$, for all $v,w\in\R^d$,
\beq
 \int |   \Phi^{-1}(V,h) -  \Phi^{-1}(W,h) | d\mu(h) \leq  K L | V - W |.
\label{jumpas31}
\eeq 
Suppose that $K\geq K L+1$.
Consider two probability densities $f^0_1,f^0_2$ on $\R^d\times\R^d$ and the corresponding solutions
$f_1,f_2$ to \eqref{eq:jump}.
It holds that for all $t\geq 0$, (here $\cT_1$ is associated to the cost function
$\varrho((x,v),(y,w))=|x-y|+|v-w|$)
$$
\cT_1 (f_1(t), f_2(t)) \leq \cT_1(f_1^0, f_2^0).
$$
\end{theorem}

\begin{proof}  As usual, we consider any probability density $F^0((x,v),(y,w))$ on  $(\R^d\times\R^d)^2$
with marginals $f^0_1$ and $f^0_2$, and we consider $F((x,v),(y,w),t)$ starting from $F^0$
and solving
\begin{eqnarray*}
&&\p_t F + v.\nabla_x F + w.\nabla_y F\\
&=& \int  \Big[F((x,\Phi(v,h)),(y,\Phi(w,h)), t)  \det (D_v \Phi(v,h))  \det (D_v \Phi(w,h)) -F((x,v),(y,w),t)\Big] d\mu(h) .
\notag
\end{eqnarray*}
This function is clearly nonnegative and has the correct marginals. For example, 
with $F_1(x,v,t) =\int F(x,y,v,w,t) dy dw$, we see that 
$$
\p_t F_1 + v.\nabla_x F_1 = \int  [F_1(x,\Phi(v,h), t)  \det (D_v \Phi(v,h)) -F_1(x,v,t)]  d\mu(h)
$$
because $\int F((x,\Phi(v,h)),(y,\Phi(w,h)), t)   \det (D_v \Phi(w,h)) dy dw = F_1(x,\Phi(v,h),t)$:
use the substitution $V=\Phi(w,h)$ (with $h$ fixed). Since $F_1(0)=f_1(0)$, we conclude that
$F_1(t)=f_1(t)$. Hence we conclude that 
$\cT_1(f_1(t),f_2(t))\leq \int\hskip-4pt\int (|x-y|+|v-w|)F((x,v),(y,w),t)dxdydvdw$.

\vip

Next, using the equation for $F$, we find with $V= \Phi(v,h)$ and $W=\Phi(w,h)$, 
\begin{eqnarray*}
&&\dis \frac d {dt} \intt (|x-y|+|v-w|) F(x,y,v,w,t)dx dy dv dw \\
&=& \dis\intt \frac{x-y}{|x-y|}\cdot(v-w)F(x,y,v,w,t) dx dy dv dw\\
&&- K \dis\intt (|x-y|+|v-w|) F(x,y,v,w,t) dx dy dv dw\\
&&+  \dis \intt \hskip-8pt \int \Big(|x-y|+| \Phi^{-1}(V,h) - \Phi^{-1}(W,h)| \Big) F((x,V),(y,W),t) 
dx dy dV dW d\mu(h)\\
&\leq& (1-K+KL) \intt |v-w| F(x,y,v,w,t)dx dy dv dw.
\end{eqnarray*}
by \eqref{jumpas31}. Since now $K\geq 1+KL$ by assumption, we deduce that
$$
\cT_1(f_1(t),f_2(t))\leq \intt (|x-y|+|v-w|)F^0((x,v),(y,w))dxdydvdw
$$
and complete the proof as usual, minimizing in $F^0$.
\end{proof} 

\begin{remark}
Fix $a>0$. If using the transport cost with weight $\varrho=a |x-y|+|v-w|$, the condition
$K\geq 1+KL$ is replaced by the condition $K > a+KL$.
\end{remark}

%-------------------------------------------------------------
\subsection{Fractional heat equation with variable coefficients}
\label{sec:FL}
%------------------------------------------------------------

Informally, the fractional Laplacian is a variant of the integral equation treated in Subsection~\ref{sec:jump}. 
However there is a particular interest when the coefficients depend on space, an example we borrow 
from~\cite{komatsu, Fournier_levyprocess}. Consider the parabolic equation with derivatives of order $\alpha \in (0,2)$
\beq \bepa
\p_t u(x,t) = \caL_\alpha [u], \qquad x \in \R, \, t \geq 0,
\\[10pt]
 \caL_\alpha ^* [\vp](x):=\dis \int [ \vp(x+\sg(x)h)-\vp(x)-h\sg(x) \vp'(x)] \frac{dh}{|h|^{1+\alpha } } .
\eepa 
\label{alLapl} \eeq

\begin{theorem} Assume that $\alpha \in (1,2)$ and that $\sg \in C^{1/\alpha}$ 
and consider two initial probability densities
$u^0_1$ and $u^0_2$ on $\R$ and the corresponding solutions $u_1$ and $u_2$ to \eqref{alLapl}.
For all $t\geq 0$,
$$
\cT_{\alpha-1} (u_1(t), u_2(t)) \leq \cT_{\alpha-1} (u_1^0, u_2^0).
$$ 
\end{theorem}

\begin{proof} 
We consider an initial probability density $v^0$ on $\R^2$ with marginals $u_1^0$ and $u_2^0$ and 
the solution $v$ to the problem (written in weak form): for all smooth $\varphi:\R^2\mapsto\R$,
\begin{eqnarray}
&&\!\!\!\frac d{dt}\intt \vp(x,y)v(x,y,t)dxdydt\label{FL:weak} \\
&\!\!\!\!=&\!\!\! \intt\! v(x,y,t) \!\int \!\! \Big(\vp(x+\sg(x)h, y+\sg(y)h)-\vp(x,y) - h[\sg(x) \partial_x\vp(x,y) + \sg(y)
\partial_y\vp_y(x,y) ] \Big) \frac{dh}{|h|^{1+\alpha}} dx dy\notag
\end{eqnarray}
starting from $v^0$. The solution is clearly nonnegative and one checks as usual that for each $t\geq 0$,
the marginals of $v(t)$ are $u_1(t)$ and $u_2(t)$: for example, we apply the above formula with 
$\varphi$ depending only on $x$ and deduce that $v_1(x,t)=\int v(x,y,t)dy$ solves the weak form of 
\eqref{alLapl}, whence $v_1=u_1$ since $v_1(0)=u_1(0)$. 
Consequently, we have $\cT_{\alpha-1} (u_1(t), u_2(t))\leq (\alpha-1)^{-1}
\int\hskip-4pt\int |x-y|^{\alpha-1}v(x,y,t)dxdy$, and, using the same arguments as usual, it suffices to show that
$$
\frac d{dt}\intt |x-y|^{\alpha-1}v(x,y,t)dxdydt \leq 0.
$$
This follows from the fact that for all $x,y\in \R$, setting $u = \frac{\sg(x)- \sg(y)}{x-y}$,
$$\bea
& \int_\R \big[ |x+\sg(x)h-y- \sg(y)h|^{\alpha-1} -|x-y|^{\alpha-1} - (\alpha -1) h[\sg(x) - \sg(y)] |x-y|^{\alpha-3}(x-y) \big] \frac{dh}{|h|^{1+\alpha}}
\\[10pt]
& \quad =  |x - y|^{\alpha-1} \int_\R \big[ |1+hu |^{\alpha-1} - 1  - (\alpha -1) h u  \big] \frac{dh}{|h|^{1+\alpha}} 
\\[10pt]
& \quad = |x - y|^{\alpha-1} |u|^\alpha \int_\R \big[ |1+h |^{\alpha-1} - 1  - (\alpha -1) h   \big] \frac{dh}{|h|^{1+\alpha}} = 0.
\eea$$
The proof of this last equality can be found in \cite[Lemma 9-(ii)]{Fournier_levyprocess},
case $a_+=a_-$ and $\beta=\alpha-1$.
Observe that  
$$
|x - y|^{\alpha-1} |u|^\alpha = \frac{|\sg(x)- \sg(y)|^\alpha}{|x-y|} \leq C_\sg
$$
so that \eqref{FL:weak} makes sense with $\varphi(x,y)=|x-y|^\alpha$, 
thanks to our regularity assumption on $\sigma$.
\end{proof} 
\vip

Here again, as in Section \ref{sec:d1}, the main technical difficulty is to prove the well-posedness of \eqref{alLapl},
in particular when $\sigma$ may degenerate. This is useful to check that the solution $v$ to the coupled equation
has the correct marginals.

%------------------------------------------------------
\section{Inhomogeneous integral equations}
\label{sec:inhomogeneous}
%------------------------------------------------------

Our next purpose is to give an example on the way to take into account $x$-dependency in IPDE models,
for instance when considering a measure $\mu(x,h)$ in the scattering equation \eqref{eq:jump}. 
We exemplify this issue  with a simple equation we borrow from \cite {FL2016}.
Consider an interval $I$ of $\R$, a rate function $d\geq 0$ defined on $I$ and some probability
density $b$ on $I$. We consider the conservative equation
\begin{equation}\label{IIE}
\p_t u(x,t) + d(x) u= b(x) A(t), \qquad A(t) = \dis \int_I d(x) u(x,t) dx
\end{equation}
starting from an initial probability density $u^0$ on $I$. We notice at once that this equation
makes sense for probability measures $u(dx,t)$ (for each $t\geq 0$, $u(dx,t)$ is a probability measure on $I$)
in the following weak sense: for all smooth $\vp:I\mapsto\R$,
\begin{equation}\label{IEEW}
\frac d{dt} \int \vp(x)u(dx,t)= \intt [\vp(z)-\vp(x)]b(z)d(x)u(dx,t)dz.
\end{equation}

\begin{theorem}
Consider two probability densities $u^0_1,u^0_2$ on $\R^d$ and the corresponding solutions
$u_1,u_2$ to \eqref{IEEW}. Under one of the two conditions (a) or (b) below, for all $t\geq 0$,
$$
\cT_\varrho (u_1(t), u_2(t)) \leq \cT_\varrho(u_1^0, u_2^0).
$$
(a) $I=\R_+$, $d(0)=0$, $d$ is increasing, $b=\delta_0$, 
and $\varrho(x,y)=|d^p(x)-d^p(y)|$ for some
$p\geq1$.

\noindent(b)  $I=\R_+$, $d(x)=\alpha x^p+\beta$ for some $\alpha,\beta\geq 0$ and $p\geq 1$,
with $\varrho(x,y)=|x^p-y^p|$, under the condition that 
$\beta \geq \alpha \int_0^\infty z^pb(z)dz$.
\end{theorem}

Other assumptions on $I$, $b$, $d$ are possible: 
it suffices that $\varrho$, $b$ and $d$ satisfy the dual inequality \eqref{ttt} below, which corresponds to~\eqref{eqp} for the heat equation with variable coefficients.

\vip

\begin{proof} We consider some probability density $v^0$ on $I^2$ with marginals $u_1^0$ and $u_2^0$ and
define the probability measure $v(dx,dy,t)$ as solving, for all smooth $\varphi:I^2\mapsto \R$, 
\begin{eqnarray}
\frac d{dt} \intt \vp(x,y)v(dx,dy,t)&=&\intt  \hskip-8pt \int\Big[\varphi(z,z)-\varphi(x,y)]b(z)\min(d(x),d(y))v(dx,dy,t)dz\notag\\
&&+\intt  \hskip-8pt \int\Big[\varphi(z,y)-\varphi(x,y)]b(z)(d(x)-d(y))_+ v(dx,dy,t)dz\notag\\
&&+\intt  \hskip-8pt \int\Big[\varphi(x,z)-\varphi(x,y)]b(z)(d(y)-d(x))_+v(dx,dy,t)dz.\label{abcd}
\end{eqnarray}
It holds true that $v(t)$ is a probability measure on $I^2$ for each $t\geq 0$, and that its
marginals are $u_1(t)$ and $u_2(t)$. For example, applying the coupling equation with
$\varphi$ depending only on $x$ and using that
$$
\min(d(x),d(y))+(d(x)-d(y))_+=d(x),
$$
one verifies that $\int_{y\in I} v(dx,dy,t)$ solves \eqref{IEEW}, whence  $\int_{y\in I} v(dx,dy,t)=u_1(dx,t)$
by uniqueness.
Hence for any cost function $\varrho:I^2\mapsto \R_+$, we have 
$\cT_\varrho(u_1(t),u_2(t))\leq\int\hskip-4pt\int \varrho(x,y) v(dx,dy,t)$.
Furthermore, we easily compute, using that $\varrho(z,z)=0$ for all $z\in I$,
that $b$ is a probability density, 
and that $\min(r,s)+(r-s)_++(s-r)_+=\max(r,s)$,
\begin{eqnarray*}
&&\dis \frac{d}{dt} \intt \varrho(x, y)  v(dx,dy,t) \dis + \intt \varrho(x, y) \max (d(x), d(y)) v(dx,dy,t)\\
&=& \dis \intt  \hskip-8pt \int \varrho(z, y) b(z) \big( d(x)-d(y) \big)_+ \; v(dx,dy,t)dz
 + \dis \intt   \hskip-8pt \int  \varrho(x, z)  b(z) \big( d(y)-d(x) \big)_+ \; v(dx,dy,t) dz .
\end{eqnarray*}
Therefore, using the same arguments as usual, the result will follows from
the fact that for all $x,y \in I$,
\begin{equation}\label{ttt}
\varrho(x, y)\max (d(x), d(y)) \geq \int [\varrho(z, y) b(z) \big( d(x)-d(y) \big)_+ 
+ \varrho(x, z)  b(z) \big( d(y)-d(x) \big)_+ dz .
\end{equation}

{\em (a)} Assuming that $I=\R_+$, that $d(0)=0$, that $d$ is increasing, that $\rho(x,y)=|d(x)-d(y)|^p$ for some
$p\geq 1$ and that $b=\delta_0$, we check that, when e.g. $x\geq y\geq 0$,
$$
(d^p(x)-d^p(y))d(x) \geq d(y)^p(d(x)-d(y)),
$$
which holds true since indeed, for any $s\geq t\geq 0$, $(s^p-t^p)s\geq t^p(s-t)$ because
$p\geq 1$.

\vip 

{\em(b)} Assume next that $I=\R_+$, $d(x)=\alpha x^p+\beta$ for some $\alpha,\, \beta \geq 0$ and $p\geq 1$
and that $\varrho(x,y)=|x^p-y^p|$. We have to verify that, for all 
$x\geq y \geq 0$,
$$
(x^p-y^p)(\alpha x^p+\beta)\geq (\alpha x^p-\alpha y^p)\int_0^\infty|z^p-y^p|b(z)dz.
$$
Setting $m=\int_0^\infty z^p b(z)dz$, it suffices to check that
$\alpha x^p+\beta \geq \alpha (m+y^p)$. This of course holds true if $\beta \geq \alpha m$.
\end{proof} 

\vip

Observe that the strong equation corresponding to the weak form \eqref{abcd} is nothing but
\begin{align*}
\p_t v + \max (d(x), d(y)) v =&  \dis b(x) \delta (x- y) \int \min(d(x') ,d(y'))  \; v(dx',dy',t)
\\
&+b(x) \int \big( d(x')- d(y) \big)_+ \; v(dx',y,t)  
\dis+b(y) \int \big( d(y')-  d(x) \big)_+ \; v(x,dy',t)  .
\end{align*}

%------------------------------------------------------
\section{Homogeneous Boltzmann equation}
\label{sec:Beq}
%------------------------------------------------------

In his seminal paper \cite{Tanaka_b}, Tanaka observed, using a probabilistic approach based on 
nonlinear Poisson-driven stochastic differential equations, 
that the homogeneous Boltzmann equation
for Maxwell mole-cules is non-expansive for the $2$-transport cost.
A deterministic proof was provided by Villani in \cite[Section 7.5.6]{VTOT}
and was extended to inelastic collisions  in \cite{BoCaBoltz}. A survey of results concerning 
homogeneous kinetic equations can be found in \cite {CaTo2007}.
The approach in \cite{VTOT,BoCaBoltz,CaTo2007} is to first derive a contractivity result 
for the gain operator 
by coupling, and then to conclude {\it via} Duhamel's principle. Here we follow Tanaka's original approach,
or rather we show how to write down his main arguments without introducing stochastic processes.

\vip

Let us also mention that \cite{rousset:hal-01020012} managed 
to study the corresponding dissipation in order to quantify the 
convergence to equilibrium of the solutions and, even more interesting, to prove the
convergence to equilibrium of the Kac particle system, with a rate of convergence not depending on
the number of particles.
\vip

The homogeneous Boltzmann writes
\beq \bepa
\p_t f(v,t) =Q(f):= \dis \int _{\R^3}\int_{\mathbb{S}^{2}}[ f(v',t)f(v_*',t) -f(v,t) f(v_*,t) ] B(\theta)  dv_* d\sg,
\\[10pt]
v' = \frac{1}{2} (v+v_*) + \frac{1}{2} |v-v_*| \sg, \quad v' _*= \frac{1}{2} (v+v_*) - \frac{1}{2}  |v-v_*| \sg,
\\[10pt]
\cos(\theta) = \frac{v-v_*}{|v-v_*|} \cdot \frac{v'-v'_*}{|v'-v'_*|}.
\eepa
\label{Boltzmann} \eeq 
The collision kernel  $B:(0,\pi)\mapsto \R_+$ is assumed to satisfy $\int_0^\pi B(\theta)d\theta=1$.
As is well-known, this equation writes, in weak form, for all 
mapping $\vp:\R^3 \mapsto \R$,
\begin{equation}\label{BoltzmannW}
\frac d{dt} \int_{\R^3} \phi(v)f(v,t)dv = \int_{\R^3\times\R^3} \int_{\mathbb{S}^2} \Big[ 
\phi(v')+\phi(v'_*)-\phi(v)-\phi(v_*)\Big] B(\theta) f(v,t)f(v_*,t)dvdv_*d\sigma.
\end{equation}

\begin{theorem}  \label{thm:Bcontract} Consider two initial probability densities $f^0_1,f^0_2$ on $\R^3$ with a finite moment of 
order $2$ and the corresponding solutions $f_1,f_2$ to \eqref{Boltzmann}.
Then, for all $t\geq 0$, one has
$$
\cT_2(f_1(t),f_2(t)) \leq \cT_2(f_1^0,f_2^0).
$$
\end{theorem}

\begin{proof} We divide the proof into four steps.
\paragraph{The global coupling.}
We fix a probability density $F^0$ on $(\R^3)^2$ with marginals $f^0_1$ and $f^0_2$ and build a coupled equation 
with the same principle as for scattering in Section~\ref{sec:integral}, that is the jump 
parameters are taken in common for the two variables,
in such a way that the post-collisional velocities are as close as possible.
We consider the solution $F(v,w,t)$, starting from $F^0$,
to the following coupling equation written in weak form: for all mapping $\Psi:\R^3\times\R^3 \mapsto \R$,
\begin{align}
\dis\frac{d}{dt}\intt_{(\R^3)^2} \Psi(v,w)F(v,w,t)dvdw
=& \dis\intt_{(\R^3\times\R^3)^2}\int_0^\pi\int_0^{2\pi}
\Big[ \Psi(v',w')+ \Psi(v'_*,w'_*) -  \Psi(v_*,w_*)- \Psi(v,w) \Big] \notag \\[4pt]
&\hskip2.5cm B(\theta) F(v,w,t)F(v_*,w_*,t)dvdw dv_*dw_* d\theta d\vp, \label{eq:Bcoupled}
\end{align}
where $v',v'_*,w',w'_*$ are built as follows.
We first choose two othonormal bases. We use the common vector $I=\frac{(v-v_*)\wedge(w-w_*)}{|(v-v_*)\wedge(w-w_*)|}$, and then 
$I_1, \, I_2$ are chosen so that  $(\frac{v-v_*}{|v-v_*|}, I,I_1)$ and  
$(\frac{w-w_*}{|w-w_*|}, I,I_2)$ are two direct orthonormal bases.
Second, we define the vectors $\sg$ and $\omega$ as
$$\begin{cases}
\sg = \cos( \theta)\frac{v-v_*}{|v-v_*|} + \sin(\theta)[ I \cos(\vp) + I_1 \sin(\vp) ] ,
\\
\omega = \cos( \theta) \frac{w-w_*}{|w-w_*|} + \sin(\theta)[ I \cos(\vp) + I_2 \sin(\vp) ].
\end{cases} $$
Finally, for $v,w,v_*,w_*\in \R^3$, $\theta\in (0,\pi)$ and $\vp \in (0,2\pi)$, the postcollisional 
velocities are defined as
\begin{gather*}
v' = \frac{1}{2} (v+v_*) + \frac{1}{2} |v-v_*| \sg, \quad v' _*= \frac{1}{2} (v+v_*) - 
\frac{1}{2}  |v-v_*| \sg,\\
w' = \frac{1}{2} (w+w_*) + \frac{1}{2} |w-w_*| \omega, \quad w' _*= \frac{1}{2} (w+w_*) - \frac{1}{2}  |w-w_*| \omega.
\end{gather*}
%\begin{gather*}v' = \frac{1}{2} (v+v_*) + \frac{1}{2} |v-v_*| \sg, \quad v' _*= \frac{1}{2} (v+v_*) - \frac{1}{2}  |v-v_*| \sg,\\
%w' = \frac{1}{2} (w+w_*) + \frac{1}{2} |w-w_*| \omega, \quad w' _*= \frac{1}{2} (w+w_*) - \frac{1}{2}  |w-w_*| \omega.\end{gather*}

\paragraph{Well posedness.} 
It is not difficult to show that \eqref{eq:Bcoupled} is well-posed in 
$L^\infty_{loc}([0,\infty),L^1(\R^3\times\R^3))$ and that the solution remains nonnegative for all times. Indeed,
applying \eqref{eq:Bcoupled} with the 
test function $\Psi(v,w)=-\ind{F(v,w)<0}$, using that 
\begin{align*}
&\Big[-\ind{F(v',w')<0}-\ind{F(v'_*,w'_*)<0}+\ind{F(v,w)<0}+\ind{F(v_*,w_*)<0}\Big]F(v,w)F(v_*,w_*)\\
\leq &4\Big[F_-(v,w)|F(v_*,w_*)|+ |F(v,w)|F_-(v_*,w_*)\Big]
\end{align*}
and that $\int_0^\pi B(\theta)d\theta=1$, we conclude that
$$
\frac{d}{dt}\intt_{(\R^3)^2} F_-(v,w,t)dvdw\leq 8\intt_{(\R^3)^2} F_-(v,w,t)dvdw \times \intt_{(\R^3)^2} |F(v,w,t)|dvdw,
$$
whence the conclusion by the Gronwall lemma.

\paragraph{Marginals.} 
Also, the marginals are correct,  that is
$\int_{\R^3}F(v,w,t)dw=f_1(v,t)$ and $\int_{\R^3}F(v,w,t)dv=f_2(w,t)$. For example concerning $f_1$,
we apply the weak coupling equation to some $\Psi$ depending only
on $v$ and we show that $\int_{\R^3}F(v,w,t)dw$ solves \eqref{BoltzmannW}. This follows 
from the fact that, when fixing $(v,w)$ and $(v_*,w_*)$, the expression between brackets in \eqref{eq:Bcoupled} 
only depends on $\sigma$, so that for any function $H:\mathbb{S}^2\mapsto \R$, we may write 
$$
\int_{\mathbb{S}^2}H(\sg)B(\theta)d\sg
=\int_0^\pi \int_0^{2\pi} H\Big(\cos( \theta)\frac{v-v_*}{|v-v_*|} 
+ \sin(\theta)[ I \cos(\vp) + I_1 \sin(\vp) ]\Big) B(\theta) d\vp d\theta.
$$
We conclude that
\[ \bea
\dis \intt_{(\R^3\times\R^3)^2}\int_0^\pi\int_0^{2\pi}
&\big[ \Psi(v')+ \Psi(v'_*)\big] B(\theta) F(v,w,t)F(v_*,w_*,t)dvdw dv_*dw_* d\theta d\vp
\\[10pt]
&\dis = \intt_{(\R^3\times\R^3)^2}\int_{\mathbb{S}^2}
\big[ \Psi(v')+ \Psi(v'_*)\big] B(\theta) F(v,w,t)F(v_*,w_*,t)dvdw dv_*dw_* d\sg
\\[10pt]
&\dis= \intt_{(\R^3\times\R^3)^2}\int_{\mathbb{S}^2}
\big[ \Psi(v')+ \Psi(v'_*) \big] B(\theta) f_1(v,t)f_1(v_*,t)dvdv_* d\sg .
\eea \]

\paragraph{The non expansion property.}
Consequently, it holds that 
$$
\cT_2(f_1(t),f_2(t))\leq \int\hskip-4pt\int_{(\R^3\times\R^3)^2} |v-w|^2 F(v,w,t)dvdw=:
h(t),
$$
 and it suffices, as usual, to show that $h'(t)\leq 0$.
For this, it suffices to verify that for all fixed
$v,w,v_*,w_* \in \R^3$, all $\theta \in (0,\pi)$,
$$
\Delta := \frac{1}{2\pi} \int_0^{2\pi} \delta d\varphi \leq 0,\quad \text{where} \quad \delta:= 
|v'-w'|^2+|v'_*-w'_*|^2-|v-w|^2-|v_*-w_*|^2.
$$
We introduce $x=v+v_*$, $y=w+w_*$, and
\begin{gather*}
a=|v-v_*|\sigma=\cos( \theta) (v-v_*) + \sin(\theta)[ I \cos(\varphi) + I_1 \sin(\varphi) ]|v-v_*| ,\\
b=|w-w_*|\omega=\cos( \theta) (w-w_*) + \sin(\theta)[ I \cos(\varphi) + I_2 \sin(\varphi) ]|w-w_*|.
\end{gather*}
We observe that
\begin{align*}
\delta=&\Big|\frac{(x-y)+(a-b)}2\Big|^2 + \Big|\frac{(x-y)-(a-b)}2\Big|^2-|v-w|^2-|v_*-w_*|^2\\
=&\frac12|x-y|^2+\frac12|a-b|^2-|v-w|^2-|v_*-w_*|^2.
\end{align*}
Using that $\int_0^{2\pi} \sin(\varphi)d\varphi=\int_0^{2\pi} \cos(\varphi)d\varphi
=\int_0^{2\pi} \sin(\varphi)\cos(\varphi) d\varphi=0$ and that $\int_0^{2\pi} \sin^2(\varphi)d\varphi
=\int_0^{2\pi} \sin^2(\varphi)d\varphi=\pi$, we find
\begin{align*}
\frac1{2\pi} \int_0^{2\pi} |a-b|^2 d \varphi=& \cos^2(\theta)\Big|(v-v_*)-(w-w_*)\Big|^2+\frac 12 \sin^2\theta 
\Big(|v-v_*|-|w-w_*|\Big)^2\\
&+ \frac 12 \sin^2(\theta) \Big||v-v_*|I_1-|w-w_*|I_2 \Big|^2.
\end{align*}
Using that $|v-v_*||w-w_*|I_1\cdot I_2 =(v-v_*)\cdot(w-w_*)$, we see that 
$$
\Big||v-v_*|I_1-|w-w_*|I_2 \Big|^2=\Big|(v-v_*)-(w-w_*)\Big|^2.
$$
Writing next
$$
\Big(|v-v_*|-|w-w_*|\Big)^2=\Big|(v-v_*)-(w-w_*) \Big|^2+2[(v-v_*)\cdot(w-w_*)- |v-v_*||w-w_*|],
$$
we get
\begin{align*}
\frac1{2\pi} \int_0^{2\pi} |a-b|^2 d \varphi=& \Big|(v-v_*)-(w-w_*)\Big|^2
+ \sin^2(\theta) [(v-v_*)\cdot(w-w_*)- |v-v_*||w-w_*|].
\end{align*}
All in all,
\begin{align*}
\Delta=& \frac12|x-y|^2+ \frac12\Big|(v-v_*)-(w-w_*)\Big|^2
+ \frac12\sin^2(\theta) [(v-v_*)\cdot(w-w_*)- |v-v_*||w-w_*|]\\
&-|v-w|^2-|v_*-w_*|^2.
\end{align*}
Recalling that $x-y=(v-w)+(v_*-w_*)$ and observing that $(v-v_*)-(w-w_*)=(v-w)-(v_*-w_*)$, we end with
\begin{equation}\label{tbc1}
\Delta=\frac12\sin^2(\theta) [(v-v_*)\cdot(w-w_*)- |v-v_*||w-w_*|].
\end{equation}
This last quantity is nonpositive, which completes the proof.
\end{proof} 

%------------------------------------------------------
\section{Porous media equation}
\label{sec:PME}
%------------------------------------------------------

We now consider the generalized porous media equation written, with $A:\R_+ \mapsto\R$ of class $C^2$, as
\beq
\p_t u - \dv (u \nabla  [A'(u)]) =0, \qquad x \in \R^d, \; t \geq 0.
\label{pmetransp} \eeq
Using a gradient flow approach  taking advantage of displacement convexity, as introduced in \cite {McC1997}, it was discovered by Otto~\cite{Otto_CPDE2001}, see also~\cite{CMV2006,OW2005,DaneriSavare}, that 
this equation is non-expansive for $\cT_2$, under a few conditions on $A$, including convexity. 
The method was improved by Bolley and Carrillo~\cite{BoCaPM} who make the contraction property an equivalence to displacement convexity for the energy functional. The coupling method, as defined in the introduction and used in the whole present paper, does not seem to apply directly.
However, using Brenier's map, the  argument  in  \cite{BoCaPM}   can somehow be 
presented in relation with a coupling. We present this argument, staying at an informal level.

\begin{theorem} [\cite{Otto_CPDE2001, BoCaPM}]
Consider some $C^2$ function $A:\R_+\mapsto\R$ such that $B(r)=\int_0^r w A''(w)dw \geq 0$ for all $r\geq 0$
and such that $r\mapsto r^{1/d-1}B(r)$ is non-decreasing.
Consider two probability densities $u_1^0,u_2^0$ on $\R^d$ and the corresponding solutions
$u_1,u_2$ to \eqref{pmetransp}. Then for all $t\geq 0$, one has
$$
\cT_2(u_1(t),u_2(t))\leq \cT_2(u_1^0,u_2^0).
$$
\end{theorem}

This applies to the porous media equation, i.e. with $A(u)=m^{-1}u^m$, as soon as $m\geq 1$.
The justification of the computation requires at least that $\int_{\R^d}B(u_1(x,t)) < \infty$.
See~\cite{Otto_CPDE2001, BoCaPM} for the rigorous proof, which 
uses that general weak solutions of the porous media equation can be approximated by smooth and positive solutions
of suitably regularized nonlinear diffusion equations.

\vip

\begin{proof} We consider Brenier's map \cite {Brenier_polar} for $u_1^0$ and $u_2^0$, i.e. a convex function $\Phi:\R^d \mapsto \R$ such that
$\cT_2(u_1^0,u_2^0)=\frac12\int_{\R^d} |x-\nabla\Phi(x)|^2u_1^0(x)dx$ and $\nabla\Phi \# u_1^0=u_2^0$.
We next consider the probability measure $v(dx,dy,t)$ (for each $t\ge 0$, $v(t)\in\mathcal{P}((\R^d)^2)$)
solving the coupling equation 
\[
\frac{\p v}{\p t} = \dv_x\big( v \nabla_x A'(u_1(x,t) \big) +\dv_y \big( v \nabla_y A'(u_2(y,t) \big)  
\]
and starting from $v^0 \in K(u^0_1,u^0_2)$ defined by the formula 
$v^0(\mathcal{A})=\int_{R^d} {\bf 1}_{\{(x,\nabla\Phi(x))\in \mathcal{A}\}} u_1^0(x)dx$
for all Borel set $\mathcal{A}\subset \R^d$, in short notation $v^0(x,y)=u_1^0(x)\delta_{\nabla\Phi(x)}(y)$. Because of this specific initial data, the equation on $v$ will not provide a global coupling as in the equations treated so far.
Again, we only use the weak form, for all smooth $\vp:\R^d\times\R^d \mapsto \R$,
$$
\frac{d}{dt}\intt\vp(x,y)v(dx,dy,t)=
-\intt \Big[\nabla_x\varphi(x,y) \cdot \nabla_x[A'(u_1(x,t))]
+\nabla_y\varphi(x,y) \cdot \nabla_y[A'(u_2(y,t))]\Big] v(dx,dy,t).
$$

\vip
One easily verifies, as usual, that for each $t\geq 0$, $v(t)$ has $u_1(t)$ and $u_2(t)$ for marginals:
for example, 
applying the weak equation of $v$ to some $\varphi$ depending only on $x$ shows that $\int_{y \in \R^d} v(dx,dy,t)$
is a (weak) solution to \eqref{pmetransp} and since it starts from $u_1^0$, we conclude by uniqueness.
As a conclusion, for all $t\geq 0$, $\cT_2(u_1(t),u_2(t)) \leq I(t)$, where
$I(t)=\frac 12 \int\hskip-4pt\int |x-y|^2 v(dx,dy,t)$.
Next,  following~\cite{BoCaPM}, we observe that
$$
I'(t)= - \intt  \big(\nabla_x [A'(u_1(x,t))]- \nabla_y [A'(u_2(y,t))]\Big)\cdot (x-y) v(dx,dy,t)=D_1(t)+D_2(t),
$$
where 
$$
D_1(t)\!=\!-\!\intt\!  \nabla_x [A'(u_1(x,t))]\cdot (x-y) v(dx,dy,t) \quad \hbox{and}\quad
D_2(t)\!=\!-\!\intt\!  \nabla_y [A'(u_2(y,t))]\cdot (y-x) v(dx,dy,t).
$$
In particular, by definition of $v^0$,
$$
D_1(0)=- \int \nabla_x [A'(u_1^0(x))]\cdot (x-\nabla\Phi(x)) u_1^0(x)dx=
-\int \nabla_x [B(u_1^0(x))]\cdot (x-\nabla\Phi(x))dx,
$$
where we recall that $B(r)=\int_0^r w A''(w)dw$.
Integrating by parts, we thus find
$$
D_1(0)= \int B(u_1^0(x)) [d-\Delta\Phi(x)]dx \leq d \int B(u_1^0(x)) [1- (\det (D^2\Phi(x)))^{1/d}]dx.
$$
This uses that for any convex function $\Phi:\R^d\mapsto\R$, we have
$d^{-1}\Delta\Phi(x)\geq   [\det (D^2\Phi(x))]^{1/d}$.  But 
since $\nabla\Phi \# u_1^0=u_2^0$, we have, for any $\vp:\R^d\mapsto\R$,
$$
\int \vp(x) u_1^0(x)dx=\int \vp((\nabla \Phi)^{-1}(y))u_2^0(y)dy=\int \vp(x)u_2(\nabla\Phi(x))\det D^2\Phi(x) dx,
$$
so that $\det D^2\Phi(x)=u_1^0(x)/u_2(\nabla\Phi(x))$ (see \cite{GuMA} for an account on this Monge-Amp\`ere equation).  All in all, we have checked that
$$
D_1(0)\leq d \int B(u_1^0(x)) \Big[1- \Big(\frac{u_1^0(x)}{u_2^0(\nabla\Phi(x))}\Big)^{1/d}\Big]dx.
$$
Proceeding similarly, we see that 
$$
D_2(0) \leq d \int B(u_2^0(y)) \Big[1- \Big(\frac{u_2^0(y)}{u_1^0((\nabla\Phi)^{-1}(y))}\Big)^{1/d}\Big]dy.
$$
Performing the substitution $x=(\nabla\Phi)^{-1}(y)$, we end up with
$$
D_2(0) \leq d \int B(u_2^0(\nabla\Phi(x))) \Big[1- \Big(\frac{u_2^0(\nabla\Phi(x))}
{u_1^0(x)}\Big)^{1/d}\Big] \frac{u_1^0(x)}{u_2^0(\nabla\Phi(x))}  dx.
$$
We thus find, with the notation $y= \nabla \Phi(x)$
\begin{eqnarray*}
\frac{I'(0)}{d} &\leq&    \int \Big[B(u_1^0(x))  \Big(1- \Big( \frac  {u_1^0(x)} {v_2( y)} \Big)^{1/d} \Big) +   B(u_2^0(y) ) \Big(1- \Big( \frac {u_2^0(y)} {u_1^0(x)} \Big)^{1/d} \Big)   \frac  {u_1^0(x)}{u_2^0(y)} \Big] dx\\
&=&\int   u_1^0(x) \Big[ \frac{B(u_1^0(x))}{u_1^0(x)} \Big(1- \Big( \frac  {u_1^0(x)} {u_2^0( y)} \Big)^{1/d} \Big) + \frac{ B(u_2^0(y) ) }{u_2^0(y)}  \Big(1- \Big( \frac {u_2^0(y)} {u_1^0(x)} \Big)^{1/d} \Big) \Big] dx\\
&=&\int   u_1^0(x) \Big[ \frac{B(u_1^0(x))}{u_1^0(x)} [u_1^0(x)]^{1/d}  -  \frac{ B(u_2^0(y) ) }{u_2^0(y)}  [u_2^0(y)]^{1/d} \Big] \Big(  [u_1^0(x)]^{-1/d}  - [u_2^0(y)]^{-1/d} \Big)  dx.
\end{eqnarray*}
Since $r\mapsto r^{1/d-1}B(r)$ is non-decreasing by assumption, we conclude that $I'(0)\leq 0$.
\vip
The above considerations hold true at any time, and not only at $t=0$. 
In other words, for all $t\geq 0$, we can find a function $I_t:[t,\infty)\mapsto \R$ such that
$I_t(t)=\cT_2(u_1(t),u_2(t))$, $I_t'(t)\leq 0$ and $\cT_2(u_1(s),u_2(s))\leq I_t(s)$ 
for all $0\leq t \leq s$.One immediately concludes that for all $t\geq 0$,
$$
\limsup_{h\downarrow 0} \frac{\cT_2(u_1(t+h),u_2(t+h))-\cT_2(u_1(t),u_2(t))}{h} \leq I_t'(t)\leq 0,
$$
so that $t\mapsto \cT_2(u_1(t),u_2(t))$ is non-increasing. 
\end{proof}

%------------------------------------------------
\section{An approach by duality}
\label{sec:duality}
%------------------------------------------------
 
In order to complete the presentation, we quickly mention another possible and original approach, based on duality. 
We consider the simplest model, i.e. the heat equation in dimension $1$, 
but all the models treated in the present paper, except the porous media equation, may be treated similarly,
with more complicated discretization procedures and more involved computations.
See also Villani \cite[pages 41-43]{VDEA} for a similar approach concerning the Vlasov equation,
in the case of the Monge-Kantorovich distance $\cT_1$, of which the dual expression is 
particularly simple.

\vip

\noindent {\bf Proof of Theorem \ref{th:H} when $d=1$ for $\rho(x,y)=|x-y|^p$ with
$p\geq 1$.} We consider two solutions
$u_1,u_2$ to \eqref{heat}, starting from probability  measures $u_1^0, u_2^0$ with finite $p$-moment.  For $h>0$, we consider
the solutions $u_{1,h},u_{2,h}$, starting from $u_1^0,u_2^0$, to the discrete heat equation
$$
\partial_t u(x,t) - \frac{1}{h^2} [u(x+h,t)+u(x-h,t)- 2 u(x,t) ] =0.
$$
 It can be written in weak form 
$$
\frac d{dt} \int_{\R^d} \vp(x)u(x,t)dx = \int_{\R^d} \frac{\varphi(x+h)+\vp(x-h)-2\vp(x)}
{h^2} u(x,t)dx.
$$
It is standard that  $u_{1,h}\to u_1$
and $u_{2,h}\to u_2$ as $h\to 0$ (in the weak topology of measures for instance). We will verify that for each $h>0$, it holds
that $\cT_p(u_{1,h}(t),u_{2,h}(t))\leq \cT_p(u_1^0,u_2^0)$, for any $p\geq 1$, and
this will complete the proof.

\vip

We fix $p\geq 1$ and
introduce the set $Q_p$ of pairs $(\varphi,\psi)$ of functions from $\R^d$ to 
$\R$ such that for all $x,y\in\R^d$, $\vp(x)+\psi(y)\leq |x-y|^p$.
For any pair of probability densities
$f,g$ on $\R^d$, the transport cost can also be expressed by duality, see \cite{V_OT2009}, as
$$
\cT_p(f,g)=\frac 1p \sup_{(\vp,\psi)\in Q_p} \Big[\int \vp(x) f(x)dx + \int \psi(y) g(y)dy\Big].
$$
For $(\varphi,\psi)\in Q_p$, we set $\Delta_{\vp,\psi}(t)=\int \vp(x)u_{1,h}(x,t)dx+
\int \psi(y)u_{2,h}(y,t)dy$. Using that $(\varphi(\cdot+h),\psi(\cdot+h))$
and $(\varphi(\cdot-h),\psi(\cdot-h))$ both belong to $Q_p$, we find
$$
\frac{d}{dt} \Delta_{\vp,\psi}(t) \leq -2h^{-2} \Delta_{\vp,\psi}(t) + 2ph^{-2} \cT_p(u_{1,h}(t),u_{2,h}(t)).
$$
This implies that
\begin{eqnarray*}
e^{2h^{-2}t} \Delta_{\vp,\psi}(t)  &\leq& \Delta_{\vp,\psi}(0) + 
2ph^{-2}\int_0^t  e^{2h^{-2}s} \cT_p(u_{1,h}(s),u_{2,h}(s))ds 
\\
&\leq& p\cT_p(u_1^0,u_2^0)+ 2ph^{-2}\int_0^t  e^{2h^{-2}s} \cT_p(u_{1,h}(s),u_{2,h}(s))ds  .
\end{eqnarray*}
Taking the supremum over all pairs $(\vp,\psi)$ in $Q_p$ and dividing by $p$,
we conclude that
$$
e^{2h^{-2}t} \cT_p(u_{1,h}(t),u_{2,h}(t))
\leq \cT_p(u_1^0,u_2^0)+
2h^{-2} \int_0^t  e^{2h^{-2}s} \cT_p(u_{1,h}(s),u_{2,h}(s))ds  .
$$
By the Gronwall lemma, we conclude that $e^{2h^{-2}t} \cT_p(u_{1,h}(t),u_{2,h}(t))
\leq e^{2h^{-2}t}  \cT_p(u_1^0,u_2^0)$ as desired.
{{\hfill $\square$}

\vip

Unfortunately, we are not able to use a similar procedure {\it directly} on the (non discretized)
heat equation. Notice however Corollary~3.3 in \cite{DpMRS16} which states, for the Kantorovich 
potentials $(\vp,\psi)$,
\[
\int \nabla u_1(x) \cdot \nabla \vp(x)dx + \int \nabla u_2(y) \cdot \nabla \psi (y) dy \geq 0.
\]
This is a type of continous form of the inequality we used here, namely 
\[
\int u_1(x) [\vp(x+h)+\vp(x-h)] dx+\int u_2(y) [\psi(y+h)+\psi(y-h)] dy\leq 2 \Big[\int u_1(x) \vp(x) dx+ \int u_2(y) \psi(y) dy \Big].
\]

\appendix

\section{Another parameterization of collisions in the Boltzmann equation}

Our construction of the global coupling for the Boltzmann equation in Section~\ref{sec:Beq} uses a standard form of collisions.  It is also standard to parametrize collisions as follows:
\beq \bepa
\p_t f(v,t) =Q(f):= \dis \int _{\R^3}\int_{\mathbb{S}^{2}}[ f(v',t)f(v_*',t) -f(v,t) f(v_*,t) ] B(\pi-2\alpha)  
dv_* d\sg,
\\[10pt]
v' = v-\tau \, (v-v_*) \cdot \tau, \quad v' _*= v_* + \tau \, (v-v_*) \cdot  \tau,
\\[10pt]
\cos(\alpha) = \frac{v-v_*}{|v-v_*|} \cdot \tau=- \frac{v'-v'_*}{|v'-v'_*|} \cdot \tau.
\eepa
\label{ap:Boltzmann} \eeq 
Our purpose is to show that one can use a similar coupling to prove the decay result 
in Theorem~\ref{thm:Bcontract}.
\vip
We introduce $I=\frac{(v-v_*)\wedge(w-w_*)}{|(v-v_*)\wedge(w-w_*)|}$, and then 
$I_1, \, I_2$ are chosen so that  $(\frac{v-v_*}{|v-v_*|}, I,I_1)$ and  
$(\frac{w-w_*}{|w-w_*|}, I,I_2)$ are two direct orthonormal bases. We then set,
for $v,v_*,w,w_* \in \R^3$, $\alpha \in (0,\pi)$ and $\varphi \in [0,2\pi)$,
$$
\begin{cases}
\tau = \cos( \alpha)\frac{v-v_*}{|v-v_*|} + \sin(\alpha)[ I \cos(\vp) + I_1 \sin(\vp) ] ,
\\
\nu = \cos( \alpha) \frac{w-w_*}{|w-w_*|} + \sin(\alpha)[ I \cos(\vp) + I_2 \sin(\vp) ].
\end{cases} 
$$
Finally, the post-collisional velocities are defined as 
\begin{gather*}
v'= v-\tau\, (v-v_*)\cdot \tau, \qquad v'_*= v_* +\tau\, (v-v_*)\cdot \tau, \\
w'= w-\nu\, (w - w_*)\cdot \nu, \qquad w'_*= w_ *+\nu\, (w - w_*)\cdot \nu .
\end{gather*}
It is also useful for later purposes to notice the equivalent formulas
\begin{gather*}
v'= v-\tau\, |v-v_*| \cos(\alpha), \qquad v'_*= v_* +\tau\,  |v-v_*| \cos(\alpha),
\\
w'= w-\nu\,  |w-w_*| \cos(\alpha), \qquad w'_*= w_ *+\nu\, |w - w_*|  \cos(\alpha).
\end{gather*}
As in Section~\ref{sec:Beq}, it suffices to verify that for all fixed
$v,w,v_*,w_* \in \R^3$, all $\alpha \in (0,\pi)$, 
$$
\Delta := \frac{1}{2\pi} \int_0^{2\pi} \big[|v'-w'|^2+|v'_*-w'_*|^2-|v-w|^2-|v_*-w_*|^2 \big] d\vp \leq 0.
$$
This computation is technical but simple. We first notice that
\begin{align*}
\tau \cdot  \nu= &\cos^2( \alpha) \frac{v-v_*}{|v-v_*|} \cdot \frac{w-w_*}{|w-w_*|} + \cos( \alpha) \sin(\vp) \frac{v-v_*}{|v-v_*|} \cdot I_2
+ \cos(\alpha) \sin(\vp)\frac{w-w_*}{|w-w_*|}  \cdot I_1 
\\
&+  \sin^2(\alpha) \cos^2(\vp) + \sin^2(\alpha) \sin^2(\vp) I_1\cdot I_2.
\end{align*}
Since $I_1\cdot I_2 =  \frac{v-v_*}{|v-v_*|} \cdot \frac{w-w_*}{|w-w_*|} $, we obtain
\begin{equation}\label{abc}
\frac{1}{2\pi} \int_0^{2\pi} \tau \cdot  \nu d\vp= \frac{1+\cos^2( \alpha)}{2} \frac{v-v_*}{|v-v_*|}  \cdot\frac{w-w_*}{|w-w_*|} + \frac{\sin^2(\alpha)}{2}  .
\end{equation}
Next, we compute
\begin{align*}
|v'-w'|^2 =& \big| v- w  - \tau \, |v-v_*| \cos(\alpha) +\nu\, |w - w_*| \cos(\alpha)  \big|^2
\\
=&  | v- w |^2 - 2(v-w)\cdot \tau \, |v-v_*| \cos(\alpha) + 2  (v - w)\cdot \nu \, |w-w_*| \cos(\alpha)
\\
&+ \big( |v-v_*| \cos(\alpha) \big)^2+ \big( |w-w_*| \cos(\alpha) \big)^2 - 2 \tau \cdot \nu |v-v_*||w-w_*|  
\cos^2(\alpha)
\end{align*}
and
\begin{align*}
|v_*'-w_*'|^2 =& \big| v_*- w_*  + \tau \, |v-v_*| \cos(\alpha) - \nu\, |w - w_*| \cos(\alpha)  \big|^2
\\
 =&  | v_*- w_* |^2 + 2(v_*-w_*)\cdot \tau \, |v-v_*| \cos(\alpha) - 2  (v_* - w_*)\cdot \nu \, |w-w_*| \cos(\alpha)
\\
& + \big( |v-v_*| \cos(\alpha) \big)^2+ \big( |w-w_*| \cos(\alpha) \big)^2 - 2 \tau \cdot \nu |v-v_*||w-w_*|  \cos^2(\alpha).
\end{align*}
Consequently, we have
\begin{align*}
&|v'-w'|^2+|v_*'-w_*'|^2 - | v- w |^2-| v_*- w_* |^2\\
=&-2 \cos(\alpha) (v-v_*-w+w_*)\cdot \tau |v-v_*|  + 2 \cos(\alpha) (v-v_*-w+w_*)\cdot \nu |w-w_*|  
\\
&+ 2 \cos^2(\alpha) [ |v-v_*|^2+|w-w_*|^2] -4 \cos^2(\alpha)  |v-v_*||w-w_*| \tau\cdot \nu.
\end{align*}
Using \eqref{abc} and that $\frac1{2\pi}\int_0^{2\pi} \tau d\varphi=\cos(\alpha)\frac{v-v_*}{|v-v_*|}$ and
$\frac1{2\pi}\int_0^{2\pi} \nu d\varphi=\cos(\alpha)\frac{w-w_*}{|w-w_*|}$, we get
\begin{align*}
\Delta=&-2 \cos^2(\alpha) (v-v_*-w+w_*)\cdot (v-v_*)  + 2 \cos^2(\alpha) (v-v_*-w+w_*)\cdot (w-w_*)  
\\
&+ 2 \cos^2(\alpha) [ |v-v_*|^2+|w-w_*|^2]\\
& -2\cos^2(\alpha)(1+ \cos^2(\alpha))(v-v_*)\cdot(w-w_*)-2\cos^2(\alpha)\sin^2(\alpha)  |v-v_*||w-w_*|.
\end{align*}
This is also
\begin{align}
\Delta=&-2 \cos^2(\alpha)\Big[ |v-v_*-w+w_*|^2 - |v-v_*|^2-|w-w_*|^2+(1+ \cos^2(\alpha))(v-v_*)\cdot(w-w_*)\\
&\hskip4cm+\sin^2(\alpha)  |v-v_*||w-w_*| \Big]\notag\\
=&-2 \cos^2(\alpha)\sin^2(\alpha)\Big[ |v-v_*||w-w_*|-(v-v_*)\cdot(w-w_*)\Big]. \label{tbc2}
\end{align}
This last quantity is nonpositive, which completes the proof. Observe that
for $\theta=\pi-2\alpha$, we have $\frac12 \sin^2(\theta)=2 \cos^2(\alpha)\sin^2(\alpha)$, so that
\eqref{tbc2} is in accordance with \eqref{tbc1}.

%
%%%%%%%%%%%%%%%%%%%%%%%%%%%%%%%%%%%
%
%%%%%% BIBLIO %%%%%%%%%%%%%%%%%%%%%%
%
%%%%%%%%%%%%%%%%%%%%%%%%%%%%%%%%%%%%
%\pagestyle{myheadings}

\end{document}